\documentclass[11pt]{article}
\usepackage[numbers,square]{natbib}
\usepackage{graphicx}
\usepackage{subcaption}
\usepackage{amsmath,amsthm,amsfonts,amssymb,commath}
\usepackage{mathtools}
\usepackage{bm}
\usepackage{booktabs}
\usepackage{xspace}
\usepackage{color}
\usepackage{stmaryrd}
\usepackage{url}
\usepackage{graphbox}
\usepackage{footnote}
\makesavenoteenv{tabular}
\makesavenoteenv{table}
\usepackage[hidelinks]{hyperref}
\hypersetup{colorlinks = true, urlcolor = blue,
  linkcolor = blue, citecolor = blue}
\usepackage{cleveref}

\usepackage[margin=0.7in]{geometry}
\makeatletter
\newcommand{\tnorm}{\@ifstar\@tnorms\@tnorm}
\newcommand{\@tnorms}[1]{%
  \left|\mkern-1.5mu\left|\mkern-1.5mu\left|
   #1
  \right|\mkern-1.5mu\right|\mkern-1.5mu\right|
}
\newcommand{\@tnorm}[2][]{%
  \mathopen{#1|\mkern-1.5mu#1|\mkern-1.5mu#1|}
  #2
  \mathclose{#1|\mkern-1.5mu#1|\mkern-1.5mu#1|}
}
\makeatother

\newtheorem{theorem}{Theorem}
\newtheorem{lemma}{Lemma}

\newtheorem{remark}{Remark}

\title{Preconditioning for the time-dependent Stokes equations:
  stage-weighted Radau IIA and weighted discontinuous Galerkin
  time-stepping}
\author{
  E. Henr\'iquez\thanks{Department of Applied Mathematics, University of
    Waterloo, ON, Canada (\url{ehenriqu@uwaterloo.ca}),
    \url{http://orcid.org/0000-0002-0243-0368}}
  \and
  J. J. Lee\thanks{Department of Mathematics, Baylor University,
    TX, USA (\url{jeonghun_lee@baylor.edu}),
    \url{https://orcid.org/0000-0001-5201-8526}}
  \and
  S. Rhebergen\thanks{Department of Applied Mathematics, University of
    Waterloo, ON, Canada (\url{srheberg@uwaterloo.ca}),
    \url{http://orcid.org/0000-0001-6036-0356}}}
\begin{document}
\maketitle
\begin{abstract}
  We determine field-of-values equivalent preconditioners for the
  time-dependent Stokes equations discretized in time by
  stage-weighted Radau IIA Runge--Kutta time-stepping and
  discontinuous Galerkin time-stepping methods. By weighing each stage
  of the Radau IIA Runge--Kutta method differently we obtain
  preconditioners that are field-of-values equivalent, with constants
  independent of viscosity, timestep, and mesh size, to the system
  matrix. Similarly, by using weighted inner products in time to
  define the discontinuous Galerkin time-stepping method, we again
  obtain field-of-values equivalent preconditioners for the resulting
  system matrix.  Numerical experiments are performed in two and
  three-dimensions which demonstrate robustness of the proposed
  preconditioners.  
\end{abstract}
\section{Introduction}
\label{s:introduction}

We consider Field-Of-Values (FOV) equivalent preconditioners for the
time-dependent Stokes equations that have been discretized in time by
stage-weighted Radau IIA Runge--Kutta (RK) and weighted discontinuous
Galerkin (DG) methods. The time-dependent Stokes problem for velocity
$u : \Omega \times I \to \mathbb{R}^d$ and pressure
$p : \Omega \times I \to \mathbb{R}$, in which
$\Omega \subset \mathbb{R}^d$, $d \in \cbr[0]{2,3}$ is a bounded
polygonal domain and $I = [0,T]$ the time-interval of interest, is
given by
\begin{equation}
  \label{eq:stokes}
  \begin{aligned}
    \partial_tu - \nabla \cdot (\nu \nabla u) + \nabla p &= f \text{ for } (x,t) \in \Omega \times I,
    &
    u &= 0 \text{ for } (x,t) \in \partial \Omega \times I,
    \\
    \nabla \cdot u &= 0 \text{ for } (x,t) \in \Omega \times I,
    &
    u(x,0) &= u_0(0) \text{ for } x \in \Omega,    
  \end{aligned}
\end{equation}
in which $\nu > 0$ is a given constant (kinematic) viscosity,
$f : \Omega \times I \to \mathbb{R}^d$ a body force, and
$u_0 : \Omega \to \mathbb{R}^d$ a given divergence-free initial
velocity field.

Preconditioners for the time-dependent Stokes equations discretized by
backward Euler have been studied in various works
\cite{gross2016robust,olshanskii2012multigrid,pearson2015block,vassilevski2013block}. Of
particular interest is the work of Cahouet and Chabard
\cite{cahouet1988some} who introduce a simple approximation to the
pressure Schur complement of the block system resulting from
discretizing \cref{eq:stokes}. This approximation is a sum of a mass
matrix scaled by $\nu^{-1}$ and a Laplacian scaled by the discrete
time step. When backward Euler time stepping is used, the discrete
form of \cref{eq:stokes} is symmetric. This symmetry is used, for
example, in \cite{mardal2004uniform,mardal2011preconditioning} to
prove that the preconditioner of Cahouet and Chabard
\cite{cahouet1988some} is parameter-robust, i.e., robust with respect
to mesh size, time step, and viscosity. They do this by finding a
uniform bound for the condition number of the preconditioned
system. See also \cite{bramble1997iterative,olshanskii2006uniform}.

When discretizing a time-dependent partial differential equation (PDE)
by a higher order time stepping method, the resulting discretization
is in general not symmetric. This complicates constructing
parameter-robust preconditioners for these methods because a uniform
condition number no longer guarantees robust convergence of Krylov
subspace methods. Instead, FOV analysis is often used to describe the
asymptotic convergence of GMRES. A functional analytic approach to
preconditioning for nonsymmetric problems is briefly discussed in
\cite{kirby2010functional}. For nonsymmetric saddle point problems we
mention the work of Loghin and Wathen \cite{loghin2004analysis} in
which they develop robust nonsymmetric FOV equivalent
preconditioners. We will use a slight modification of this framework
to determine robust FOV equivalent preconditioners for system matrices
resulting from Radau IIA RK \cite{EhleThesis,hairer1999stiff} and
higher-order DG \cite{jamet1978galerkin,thomee2007galerkin} time
stepping methods applied to \cref{eq:stokes}. Our preconditioners are
higher-order variations of the Cahouet and Chabard
\cite{cahouet1988some} preconditioner. Let us also mention
\cite{clines2022efficient,leveque2026augmented,mardal2007order,pazner2017stage,southworth2022fastb,southworth2022fasta}
for alternative strategies that have been proposed to precondition
implicit RK methods applied to time-dependent PDEs.

In this paper we first consider preconditioners for Radau IIA RK time
stepping methods. By weighing each stage of the RK scheme by a
suitably chosen parameter we are able to show that the discrete weak
form of $\partial_tu - \nabla \cdot (\nu \nabla u)$ in
\cref{eq:stokes} when using an $s$-stage Radau IIA RK discretization
is coercive. This allows us to design a simple and robust block
preconditioner (one block per stage) for the system matrix.

We next consider preconditioners for higher-order \textit{weighted} DG
time stepping methods. The idea of using DG time stepping with
weighted inner products in time to discretize time-dependent PDEs has
previously been used, for example, in
\cite{french1993space,sosa2020space}. Weighted inner products in time
have also been used to analyze higher-order DG time stepping methods
\cite{chrysafinos2010discontinuous,kirk2023analysis,wang2026space}. We
use weighted inner products in time to discretize \cref{eq:stokes} in
time and show that the discrete weak form of
$\partial_tu - \nabla \cdot (\nu \nabla u)$ in \cref{eq:stokes} is
coercive. This result is in contrast to using traditional DG
time-stepping methods in which the discrete weak form of
$\partial_tu - \nabla \cdot (\nu \nabla u)$ is inf-sup stable but not
necessarily coercive. We will use the coercivity of this weak form to
design a robust preconditioner for the weighted DG time stepping
method.

Traditional (non-weighted) DG time stepping methods are equivalent to
Radau IIA RK time stepping methods
\cite{Akrivis:2022,Akrivis:2011}. We show that the weighted DG time
stepping method results in a new RK method, one in which the RK
time-discretization of $\partial_tu - \nabla \cdot (\nu \nabla u)$ in
\cref{eq:stokes} is coercive. This new RK method provides hints as to
how to choose the weights in the weighted Radau IIA RK time stepping
method.

We organize this paper as follows. In \cref{s:saddlepointproblem} we
summarize, and slightly modify, the framework presented in
\cite{loghin2004analysis} for FOV-equivalent preconditioners for
nonsymmetric saddle-point problems. In \cref{s:t-stokes} we present
and analyze the stage-weighted Radau IIA RK, weighted DG, and new RK
time-discretizations of \cref{eq:stokes}. Parameter-robust FOV
equivalent preconditioners are presented in
\cref{s:preconditioning}. Numerical examples in two and three
dimensions in \cref{s:numex} verify the robustness of the proposed
preconditioners.

\section{Preconditioning nonsymmetric saddle-point problems}
\label{s:saddlepointproblem}

Before determining a preconditioner for discretizations of
\cref{eq:stokes}, we analyze FOV-equivalent preconditioners for a
general nonsymmetric saddle-point problem. Let $\cbr[0]{V_h}$ and
$\cbr[0]{Q_h}$ be families of finite element spaces corresponding to
shape regular meshes of size $h$. Let $V_h$ be endowed with an inner
product $(\cdot, \cdot)_{V_h}$ and denote by $\norm[0]{\cdot}_{V_h}$
the norm induced by this inner product. Similarly, let $Q_h$ be
endowed with an inner product $(\cdot, \cdot)_{Q_h}$ and denote its
induced norm by $\norm[0]{\cdot}_{Q_h}$. We consider the problem: Find
$(u_h, p_h) \in V_h \times Q_h$, such that
\begin{subequations}
  \label{eq:saddlepoint}
  \begin{align}
    a(u_h,v_h) + b(v_h,p_h) &= \mathcal{F}(v_h) && \forall v_h \in V_h,
    \\
    -b(u_h,q_h) &= \mathcal{G}(q_h) && \forall q_h \in Q_h,
  \end{align}
\end{subequations}
where $\mathcal{F} : V_h \to \mathbb{R}$ and
$\mathcal{G} : Q_h \to \mathbb{R}$ are bounded linear functionals. We
assume there exist positive constants $c_a$ and $\alpha$ such that the
bilinear form $a: V_h \times V_h \to \mathbb{R}$, which need not be
symmetric, satisfies
\begin{subequations}
  \label{eq:aconditions}
  \begin{align}
    \label{eq:a-boundedness}
    a(u_h,v_h) &\le c_a \norm[0]{u_h}_{V_h} \norm[0]{v_h}_{V_h} && \forall u_h,v_h \in V_h,
    \\
    \label{eq:abs-coercivity}
    a(v_h, v_h) &\ge \alpha \norm[0]{v_h }_{V_h}^2 && \forall v_h \in V_h,
  \end{align}
\end{subequations}
and that there exist positive constants $c_b$ and $\beta$ such that
the bilinear form $b: V_h \times Q_h \to \mathbb{R}$ satisfies
\begin{subequations}
  \label{eq:bconditions}
  \begin{align}
    \label{eq:b-boundedness}
    &b(v_h,q_h) \le c_b \norm[0]{v_h}_{V_h} \norm[0]{q_h}_{Q_h}
    \qquad \forall v_h \in V_h, \ q_h \in Q_h,
    \\
    \label{eq:abs-infsup}
    &\inf_{q_h \in Q_h} \sup_{v_h \in V_h}
    \frac{b(v_h, q_h)}{\norm[0]{v_h}_{V_h} \norm[0]{q_h}_{Q_h}}
    \ge \beta.
  \end{align}
\end{subequations}
The assumptions in \cref{eq:bconditions,eq:aconditions} are sufficient
for the saddle-point problem \cref{eq:saddlepoint} to have a unique
solution \cite[Theorem 4.2.1]{boffi2013mixed}.

Let $\langle {\bf x}, {\bf y} \rangle = {\bf x}^T {\bf y}$
$\forall {\bf x},{\bf y} \in \mathbb{R}^n$ be the Euclidean inner
product on $\mathbb{R}^n$. Let $M \in \mathbb{R}^{n \times n}$ be any
symmetric and positive definite matrix. We will write
$\langle {\bf x}, {\bf y} \rangle_M := \langle M{\bf x}, {\bf y}
\rangle$ $\forall {\bf x},{\bf y} \in \mathbb{R}^n$. This inner
product induces the norm
$\norm[0]{{\bf x}}_M = \langle {\bf x}, {\bf x}
\rangle_M^{1/2}$. Furthermore, we say that the matrices
$A,B \in \mathbb{R}^{n \times n}$ are ($M$-)FOV equivalent if there
exist positive constants $c_{\star}$ and $\bar{c}_{\star}$ such that
\begin{equation*}
  c_{\star}
  \le \inf_{{\bf 0} \ne {\bf x} \in \mathbb{R}^n} \frac{\langle AB^{-1}{\bf x}, {\bf x}\rangle_M}{\langle {\bf x}, {\bf x}\rangle_M}, \qquad
  \bar{c}_{\star}
  \le \inf_{{\bf 0} \ne {\bf x} \in \mathbb{R}^n} \frac{\langle BA^{-1}{\bf x}, {\bf x} \rangle_M}{\langle{\bf x}, {\bf x}\rangle_M}.
\end{equation*}

Let $n_v = \dim V_h$, $n_q = \dim Q_h$, $n = n_v + n_q$,
$\cbr[0]{\phi_i}_{i=1}^{n_v}$ be a basis for $V_h$ and
$\cbr[0]{\psi_i}_{i=1}^{n_q}$ be a basis for $Q_h$. Let $v_h \in V_h$
such that $v_h = \sum_{i=1}^{n_v} {\bf v}_i \phi_i$ and define
${\bf v} \in \mathbb{R}^{n_v}$ to be the resulting vector of expansion
coefficients. Likewise, let $q_h \in Q_h$ such that
$q_h = \sum_{i=1}^{n_q} {\bf q}_i \psi_i$ and define
${\bf q} \in \mathbb{R}^{n_q}$ to be the resulting vector of expansion
coefficients. We define the matrices $A$ and $B$ such that
$\langle A{\bf v}, {\bf v} \rangle = a_h(v_h,v_h)$ and
$\langle B{\bf v}, {\bf q}\rangle = b_h(v_h,q_h)$, and the vectors
${\bf f}$ and ${\bf g}$ such that
$\langle {\bf f}, {\bf v} \rangle = \mathcal{F}(v_h)$ and
$\langle {\bf g}, {\bf q} \rangle = \mathcal{G}(q_h)$. We may then
write \cref{eq:saddlepoint} as finding ${\bf u} \in \mathbb{R}^{n_v}$
and ${\bf p} \in \mathbb{R}^{n_q}$ such that
\begin{equation}
  \label{eq:saddlepointMatrix}
  \mathcal{A} 
  \begin{bmatrix}
    {\bf u} \\ {\bf p}
  \end{bmatrix}
  =
  \begin{bmatrix}
    A & B^T \\ -B & 0
  \end{bmatrix}
  \begin{bmatrix}
    {\bf u} \\ {\bf p}
  \end{bmatrix}
  =
  \begin{bmatrix}
    {\bf f} \\ {\bf g}
  \end{bmatrix}.
\end{equation}
Define the matrices $P_v$ and $P_q$ such that
\begin{equation*}
  \langle P_v {\bf u}, {\bf v} \rangle = (u_h, v_h)_{V_h} \quad \forall u_h,v_h \in V_h,
  \qquad
  \langle P_q {\bf p}, {\bf q} \rangle = (p_h, q_h)_{Q_h} \quad \forall p_h,q_h \in Q_h,
\end{equation*}
and 
\begin{equation}
  \label{eq:preconditionerPULW}
  \widetilde{\mathcal{P}} :=
  \begin{bmatrix}
    P_v & 0 \\ 0 & P_q
  \end{bmatrix},
  \qquad
  \widetilde{\mathcal{P}}_U :=
  \begin{bmatrix}
    P_v & B^T \\ 0 & \rho P_q
  \end{bmatrix},
\end{equation}
with $\rho$ a positive parameter. Under the assumptions
\cref{eq:aconditions,eq:bconditions}, Loghin and Wathen
\cite{loghin2004analysis} show that $\widetilde{\mathcal{P}}_U$ is
($\widetilde{\mathcal{P}}^{-1}$)-FOV equivalent to $\mathcal{A}$,
provided $\rho \le \rho_0$ where $\rho_0$ is a constant such that
$\norm[0]{AP_v^{-1}-I}_{P_v^{-1}} \le \rho_0$, see \cite[Theorem 3.9,
Remark 3.3]{loghin2004analysis}. A consequence of this result is that
the residuals of the preconditioned system, using preconditioner
$\widetilde{\mathcal{P}}_U$, at the $k$-th GMRES iteration,
${\bf r}^k$, satisfy \cite{starke1997field},
\begin{equation}
  \label{eq:gmresconv}
  \frac{\norm[0]{{\bf r}^k}_{\widetilde{\mathcal{P}}^{-1}}}{\norm[0]{{\bf r}^0}_{\widetilde{\mathcal{P}}^{-1}}} \le
  (1 - c_{\star} \bar{c}_{\star})^{k/2}.
\end{equation}

\begin{remark}
  \label{rem:LogWathgeneral}
  Loghin and Wathen \cite{loghin2004analysis} also prove FOV
  equivalence for the more general case in which
  \cref{eq:abs-coercivity} is replaced by $a(\cdot, \cdot)$ being
  inf-sup stable. In this case, they show that
  $\widetilde{\mathcal{P}}_U$, but with $P_v$ replaced by $A$ or a
  matrix that is ($P_v^{-1}$)-FOV equivalent to $A$, is also
  ($\widetilde{\mathcal{P}}^{-1}$)-FOV equivalent to $\mathcal{A}$,
  see \cite[Theorem 3.7 and 3.9]{loghin2004analysis}.
\end{remark}

In this work we consider a slight modification of
\cref{eq:preconditionerPULW}. Let
\begin{equation}
  \label{eq:preconditionerPU}
  \mathcal{P} :=
  \begin{bmatrix}
    \omega P_v & 0 \\ 0 & \rho P_q
  \end{bmatrix},
  \qquad
  \mathcal{P}_U :=
  \begin{bmatrix}
    \omega P_v & B^T \\ 0 & \rho P_q
  \end{bmatrix}.
\end{equation}
For completeness, we redo the proofs to show that $\mathcal{P}_U$ is
($\mathcal{P}^{-1}$-)FOV-equivalent to $\mathcal{A}$ for $\omega$ and
$\rho$ sufficiently large positive parameters. Our proofs are
different from \cite{loghin2004analysis}; the assumption that
$a(\cdot, \cdot)$ is coercive \cref{eq:abs-coercivity} simplifies the
proofs. We have the following result.

\begin{theorem}
  \label{thm:fovequiv}
  Let $c_a,c_b,\alpha,\beta$ be the constants from
  \cref{eq:bconditions,eq:aconditions} and let $\omega$ and $\rho$ be
  such that $\rho \omega > \tfrac{1}{4}c_a^2\alpha^{-1}$ and
  $\omega > \tfrac{1}{2}c_a^2\alpha^{-1}$. Then, the matrices
  $\mathcal{A}$ and $\mathcal{P}_U$, as defined in
  \cref{eq:saddlepointMatrix,eq:preconditionerPU} are
  ($\mathcal{P}^{-1}$-)FOV-equivalent with
  \begin{equation*}
    c_{\star} = \frac{2\sqrt{\omega \alpha} - c_a/\sqrt{\rho}}{2\sqrt{\omega^3 \alpha^{-1}} \max\{1, \alpha^{-2}c_a^3\beta^{-2} \}},
    \qquad
    \bar{c}_{\star} = \frac{\min\{2\alpha\omega - c_a^2, 1\}}{2\max\{2 c_a^2+ \omega \rho^{-1} c_b^2, 2\}}. 
  \end{equation*}
\end{theorem}
\begin{proof}
  See Appendix \ref{ap:proofFOVequiv}.
\end{proof}

\begin{remark}
  \label{rem:uniformwellposed}
  Assume that the inner products $(\cdot, \cdot)_{V_h}$ and
  $(\cdot, \cdot)_{Q_h}$ are defined such that the constants
  $\alpha, \beta, c_a, c_b$ in \cref{eq:aconditions,eq:bconditions}
  are independent of mesh size and other problem parameters, i.e., are
  uniform constants. Then \cref{thm:fovequiv} and \cref{eq:gmresconv}
  imply that asymptotic convergence of GMRES of the preconditioned
  system is independent of mesh size and other problem parameters.
\end{remark}

\section{The fully-discrete time-dependent Stokes problem}
\label{s:t-stokes}

To discretize \cref{eq:stokes} in time we consider both Radau IIA RK
time-stepping and DG time-stepping methods. For simplicity, we
consider a conforming and inf-sup stable finite element spatial
discretization, i.e., $V_h \subset [H_0^1(\Omega)]^d$ and
$Q_h \subset L_0^2(\Omega)$. In particular, we will use the
Taylor--Hood element. Let $(u,v) := \int_{\Omega}uv \dif x$ be the
$L^2$-inner product on $\Omega$, $a(u,v) := \nu(\nabla u, \nabla v)$,
and $b(v,q) := - (\nabla \cdot v, q)$. To analyze the discretizations,
we consider the following inner products on $V_h$ and $Q_h$,
respectively:
\begin{equation}
  \label{eq:innerprodvq}
  \begin{split}
    (u_h, v_h)_v
    &:= \tau (u_h, v_h) + \nu (\nabla u_h, \nabla v_h),
    \\
    (p_h, q_h)_q
    &:= \inf_{\substack{p_h = p_{h,1} + p_{h,0} \\ q_h = q_{h,1} + q_{h,0}}}
    \del[1]{\tau^{-1}(\nabla p_{h,1}, \nabla q_{h,1}) + \nu^{-1}(p_{h,0}, q_{h,0})}.    
  \end{split}
\end{equation}
The norms induced by these inner products are denoted by
$\norm[0]{\cdot}_v$ and $\norm[0]{\cdot}_q$. As discussed in
\cite{mardal2004uniform}, $\norm[0]{\cdot}_v$ is a norm on
$\tau [L^2(\Omega)]^d \cap \nu [H_0^1(\Omega)]^d$ and
$\norm[0]{\cdot}_q$ is a norm on
$\tau^{-1}(H^1(\Omega) \cap L_0^2(\Omega)) +
\nu^{-1}L_0^2(\Omega)$. It is shown in \cite[Section
5.1]{mardal2004uniform}, that when using the Taylor--Hood element,
there exists a uniform constant $\beta > 0$ such that
\begin{equation}
  \label{eq:infsupmardal}
  \sup_{v_h \in V_h} \frac{b(v_h, q_h)}{\norm[0]{v_h}_{v}} \ge \beta \norm[0]{q_h}_{q} \qquad
  \forall q_h \in Q_h.
\end{equation}

The triangulation of the domain $\Omega$ is denoted by
$\mathcal{T}_h := \cbr[0]{K}$, which consists of simplices $K$. We
partition the time-interval into $N$ subintervals of length
$\Delta t = T/N$ and define $t_n := n \Delta t$ for $n=0,\hdots,N$,
$I_n = [t_n,t_{n+1}]$ for $n=0,\hdots,N-1$, and $\tau := 1/\Delta t$.

\subsection{Radau IIA Runge--Kutta time-stepping}
\label{ss:radauiia}

We first discretize \cref{eq:stokes} in time by a Radau IIA RK
time-stepping method \cite{EhleThesis,hairer1999stiff}. Let
$\mathbb{A}_r = [a_{ij}^r] \in \mathbb{R}^{s \times s}$ be the Radau
IIA RK matrix, $0 < c_1^r < \hdots < c_s^r = 1$ the Radau IIA RK
nodes, and let $t_{n,i} = t_n + c_i^r \Delta t$. Define
$f_j(x) = f(x,t_{n,i})$, denote by $(u_i,p_i)$ the approximation to
$(u,p)$ at $t_{n,i}$, and denote by $(u^n,p^n)$ the approximation to
$(u,p)$ at $t_n$. Discretizing \cref{eq:stokes} in time we have, given
$u^n$, that $(u^{n+1},u^{n+1})$ is found by solving (see, for example,
\cite{hairer2010solving}):
\begin{equation}
  \label{eq:stokes-radauiia}
  u_i = u_n + \Delta t \sum_{j=1}^s a_{ij}^r\sbr[1]{f_j + \nabla \cdot (\nu \nabla u_j) - \nabla p_j},
  \qquad
  \nabla \cdot u_i = 0,
\end{equation}
for $1 \le i \le s$ and setting $(u^{n+1},p^{n+1}) =
(u_s,p_s)$. Multiplying \cref{eq:stokes-radauiia} by the inverse of
the RK matrix and discretizing also in space, we obtain the fully
discrete problem: For $n=0,\hdots, N-1$, find
$\boldsymbol{u}_h^n := (U_1,\hdots,U_s) \in \boldsymbol{V}_h :=
[V_h]^s$ and
$\boldsymbol{p}_h^n := (P_1,\hdots,P_s) \in \boldsymbol{Q}_h :=
[Q_h]^s$ such that
\begin{equation}
  \label{eq:stokes-radauiiaTH}
  \boldsymbol{\widetilde{a}}_r(\boldsymbol{u}_h, \boldsymbol{v}_h)
  + \boldsymbol{\widetilde{b}}_r(\boldsymbol{v}_h, \boldsymbol{p}_h)
  - \boldsymbol{\widetilde{b}}_r(\boldsymbol{u}_h, \boldsymbol{q}_h)
  = \sum_{i=1}^s \sbr[3]{(f_i, V_i)
  + \tau \sum_{j=1}^s (\mathbb{A}_r^{-1})_{ij}(U_n, V_i)},
\end{equation}
for all $\boldsymbol{v}_h := (V_1,\hdots,V_s) \in \boldsymbol{V}_h$
and $\boldsymbol{q}_h := (Q_1,\hdots,Q_s) \in \boldsymbol{Q}_h$, in
which
\begin{equation*}
  \boldsymbol{\widetilde{a}}_r(\boldsymbol{u}_h, \boldsymbol{v}_h)
  := \sum_{i=1}^s \sbr[3]{\tau \sum_{j=1}^s (\mathbb{A}^{-1}_r)_{ij}(U_j, V_i)
    +  a(U_i, V_i)},
  \quad
  \boldsymbol{\widetilde{b}}_r(\boldsymbol{v}_h, \boldsymbol{p}_h)
  := \sum_{i=1}^s b(V_i, P_i).
\end{equation*}
We then set $(u_h^{n+1},p_h^{n+1}) = (U_s,P_s)$. The bilinear form
$\boldsymbol{\widetilde{a}}_r(\cdot, \cdot)$ cannot be shown to
satisfy \cref{eq:abs-coercivity}. However, it is known that the RK
matrix $\mathbb{B}$ of an A-stable irreducible RK method of which
$\mathbb{B}$ is invertible is weakly positive definite, i.e., that
there exists a matrix $\widehat{W} \in \mathbb{R}^{s \times s}$ such
that $\langle \widehat{W}{\bf x}, {\bf x}\rangle > 0$ and
$\langle \widehat{W}\mathbb{B}{\bf x}, {\bf x} \rangle > 0$ for all
${\bf x} \in \mathbb{R}^s$, see \cite[Definition 3.1, Lemma
3.2]{mardal2007order}. Assuming therefore that there exists a
$\widehat{W} := \text{diag}(\hat{w}_1,\hdots,\hat{w}_s)$ with positive
$\hat{w}_i \in \mathbb{R}$ such that
$\langle \widehat{W}\mathbb{A}_r{\bf x}, {\bf x} \rangle > 0$, we
multiply each stage in \cref{eq:stokes-radauiiaTH} by $\hat{w}_i$ and
introduce a stage-weighted Radau IIA + Taylor--Hood discretization:
For $n=0,\hdots, N-1$, find
$\boldsymbol{u}_h^n := (U_1,\hdots,U_s) \in \boldsymbol{V}_h :=
[V_h]^s$ and
$\boldsymbol{p}_h^n := (P_1,\hdots,P_s) \in \boldsymbol{Q}_h :=
[Q_h]^s$ such that
\begin{equation}
  \label{eq:stokes-wradauiia}
  \boldsymbol{a}_r(\boldsymbol{u}_h, \boldsymbol{v}_h)
  + \boldsymbol{b}_r(\boldsymbol{v}_h, \boldsymbol{p}_h)
  - \boldsymbol{b}_r(\boldsymbol{u}_h, \boldsymbol{q}_h)
  = \sum_{i=1}^s \hat{w}_i \sbr[3]{(f_i, V_i)
  + \tau \sum_{j=1}^s (\mathbb{A}_r^{-1})_{ij}(U_n, V_i)},
\end{equation}
for all $\boldsymbol{v}_h := (V_1,\hdots,V_s) \in \boldsymbol{V}_h$
and $\boldsymbol{q}_h := (Q_1,\hdots,Q_s) \in \boldsymbol{Q}_h$, in
which
\begin{equation*}
  \boldsymbol{a}_r(\boldsymbol{u}_h, \boldsymbol{v}_h)
  := \sum_{i=1}^s \hat{w}_i\sbr[3]{\tau \sum_{j=1}^s (\mathbb{A}_r^{-1})_{ij}(U_j, V_i)
    +  a(U_i, V_i)},
  \
  \boldsymbol{b}_r(\boldsymbol{v}_h, \boldsymbol{p}_h)
  := \sum_{i=1}^s \hat{w}_i b(V_i, P_i).
\end{equation*}

In what follows we will prove that $\boldsymbol{a}_r(\cdot, \cdot)$
and $\boldsymbol{b}_r(\cdot, \cdot)$ satisfy
\cref{eq:aconditions,eq:bconditions} assuming the existence of the
aforementioned matrix $\widehat{W}$. We discuss how to choose
$\widehat{W}$ in \cref{s:numex}. We introduce the following inner
products associated with $\boldsymbol{V}_h$ and $\boldsymbol{Q}_h$:
\begin{equation*}
  (\boldsymbol{u}_h, \boldsymbol{v}_h)_{r,\boldsymbol{v}}
  := \sum_{i=1}^s \hat{w}_i(u_i, v_i)_v,
  \qquad
  (\boldsymbol{p}_h, \boldsymbol{q}_h)_{r,\boldsymbol{q}}
  := \sum_{i=1}^s \hat{w}_i(p_i, q_i)_q,
\end{equation*}
respectively, with induced norms denoted by
$\norm[0]{\cdot}_{r,\boldsymbol{v}}$ and
$\norm[0]{\cdot}_{r,\boldsymbol{q}}$. Here $(\cdot, \cdot)_v$ and
$(\cdot, \cdot)_q$ are defined in \cref{eq:innerprodvq}.

\begin{lemma}[Coercivity]
  \label{lem:coercivity-wradauiia}
  Let
  $\widehat{W} := \text{diag}(\hat{w}_1,\hdots,\hat{w}_s) \in
  \mathbb{R}^{s \times s}$. If the $\hat{w}_i$ are such that
  $C_r := \tfrac{1}{2}( \widehat{W}\mathbb{A}_r^{-1} +
  (\widehat{W}\mathbb{A}_r^{-1})^T)$ is symmetric positive definite,
  then there exists a uniform constant $\alpha > 0$ such that
  \begin{equation*}
    \boldsymbol{a}_r(\boldsymbol{v}_h, \boldsymbol{v}_h) \ge \alpha \norm[0]{\boldsymbol{v}_h}_{r,\boldsymbol{v}}^2
    \qquad \forall \boldsymbol{v}_h \in \boldsymbol{V}_h.
  \end{equation*}
\end{lemma}
\begin{proof}
  Since
  $C_r = \tfrac{1}{2}( \widehat{W}\mathbb{A}_r^{-1} +
  (\widehat{W}\mathbb{A}_r^{-1})^T)$ is symmetric positive definite,
  we have by the spectral theorem that
  $[C_r]_{ij} = \sum_{m=1}^sZ_{im}\lambda_mZ_{jm}$ with
  $\lambda_m > 0$ the eigenvalues of $C_r$ and
  $Z \in \mathbb{R}^{s\times s}$ an orthogonal matrix.  Then, since
  $\sum_{i=1}^s\sum_{j=1}^s \tau \hat{w}_i (\mathbb{A}_r^{-1})_{ij}(V_j, V_i)
  = \sum_{i=1}^s\sum_{j=1}^s \tau [C_r]_{ij}(V_j, V_i)$,
  \begin{align*}
    \boldsymbol{a}_r(\boldsymbol{v}_h, \boldsymbol{v}_h)
    &= \tau \sum_{m=1}^s\lambda_m\sum_{i=1}^s \sum_{j=1}^s Z_{jm} Z_{im} (V_j, V_i)
      + \nu \sum_{i=1}^s \hat{w}_i \norm[0]{\nabla V_i}_{L^2(\Omega)}^2
    \\
    &\ge \lambda_{\min} \tau \sum_{i=1}^s \sum_{j=1}^s (\sum_{m=1}^sZ_{mj}Z_{mi}) (V_j, V_i)
      + \nu \sum_{i=1}^s \hat{w}_i \norm[0]{\nabla V_i}_{L^2(\Omega)}^2
    \\
    &= \lambda_{\min} \tau \sum_{i=1}^s \sum_{j=1}^s \delta_{ij} (V_j, V_i)
      + \nu \sum_{i=1}^s \hat{w}_i \norm[0]{\nabla V_i}_{L^2(\Omega)}^2
    \\
    &= \lambda_{\min} \tau \sum_{i=1}^s \norm[0]{V_i}_{L^2(\Omega)}^2
      + \nu \sum_{i=1}^s \hat{w}_i \norm[0]{\nabla V_i}_{L^2(\Omega)}^2
    \\
    &\ge \min\del[2]{\frac{\lambda_{\min}}{\max_{1\le k\le s}\hat{w}_k}, 1}\norm[0]{\boldsymbol{v}_h}_{r,\boldsymbol{v}}^2,
  \end{align*}
  in which $\lambda_{\min}$ is the smallest eigenvalue of $C_r$.
\end{proof}

\begin{lemma}[Stability]
  \label{lem:stability-wradauiia}
  There exists a uniform constant $\beta > 0$ such that
  \begin{equation*}
    \sup_{\boldsymbol{v}_h \in \boldsymbol{V}_h} \frac{\boldsymbol{b}_r(\boldsymbol{v}_h, \boldsymbol{q}_h)}{\norm[0]{\boldsymbol{v}_h}_{r,\boldsymbol{v}}}
    \ge \beta \norm[0]{\boldsymbol{q}_h}_{r,\boldsymbol{q}} \qquad
    \forall \boldsymbol{q}_h \in \boldsymbol{Q}_h.
  \end{equation*}
\end{lemma}
\begin{proof}
  By \cref{eq:infsupmardal}, given $Q_i \in Q_h$, there exists
  $\tilde{V}_i \in V_h$ such that
  $b(\tilde{V}_i, Q_i) = \norm[0]{Q_i}_q^2$ and
  $\norm[0]{\tilde{V}_i}_v \le \beta^{-1} \norm[0]{Q_i}_q$. Define
  $\boldsymbol{\tilde{v}}_h := (\tilde{V}_1, \hdots, \tilde{V}_s) \in
  \boldsymbol{V}_h$. We find
  \begin{align*}
    \norm[0]{\boldsymbol{\tilde{v}}_h}_{r,\boldsymbol{v}}^2
    &= \sum_{i=1}^s \hat{w}_i \norm[0]{\tilde{V}_i}_v^2
      \le \sum_{i=1}^s \hat{w}_i \beta^{-2} \norm[0]{Q_i}_q^2
      = \beta^{-2} \norm[0]{\boldsymbol{q}_h}_{r,\boldsymbol{q}}^2,
    \\
    \boldsymbol{b}_r(\boldsymbol{\tilde{v}}_h, \boldsymbol{q}_h)
    &= \sum_{i=1}^s \hat{w}_i b(\tilde{V}_i, q_i)
      = \sum_{i=1}^s \hat{w}_i \norm[0]{Q_i}_q^2
      = \norm[0]{\boldsymbol{q}_h}_{r,\boldsymbol{q}}^2.      
  \end{align*}
  The result follows.
\end{proof}

\begin{lemma}[Continuity]
  \label{lem:continuity-wradauiia}
  There exist uniform constants $c_a,c_b > 0$ such that
  \begin{subequations}
    \begin{align}
      \label{eq:boldabound-wradauiia}
      |\boldsymbol{a}_r(\boldsymbol{u}_h, \boldsymbol{v}_h)|
      &\le c_a \norm[0]{\boldsymbol{u}_h}_{r,\boldsymbol{v}} \norm[0]{\boldsymbol{v}_h}_{r,\boldsymbol{v}}
      && \forall \boldsymbol{u}_h, \boldsymbol{v}_h \in \boldsymbol{V}_h,
      \\
      \label{eq:boldbbound-wradauiia}
      |\boldsymbol{b}_r(\boldsymbol{v}_h, \boldsymbol{q}_h)|
      &\le c_b \norm[0]{\boldsymbol{v}_h}_{r,\boldsymbol{v}} \norm[0]{\boldsymbol{q}_h}_{r,\boldsymbol{q}}
      && \forall \boldsymbol{v}_h \in \boldsymbol{V}_h,\ \boldsymbol{q}_h \in \boldsymbol{Q}_h.
    \end{align}    
  \end{subequations}
\end{lemma}
\begin{proof}
  Using a triangle inequality and Cauchy--Schwarz,
  \begin{equation}
    \label{eq:aboundp-wradauiia}
    |\boldsymbol{a}_r(\boldsymbol{u}_h, \boldsymbol{v}_h)|
    \le
    \envert[3]{ \sum_{i=1}^s \sum_{j=1}^s \tau \hat{w}_i (\mathbb{A}_r^{-1})_{ij}(U_j, V_i) }
    + \envert[3]{ \sum_{i=1}^s \hat{w}_i \nu (\nabla U_i, \nabla V_i) }
    =: T_1 + T_2.
  \end{equation}
  We bound $T_1$ and $T_2$ separately. First, using Cauchy--Schwarz,
  \begin{equation*}
    \begin{split}
      T_1
      &\le \max_{1 \le \ell,k \le s}[(\mathbb{A}^{-1}_r)_{\ell k}] \max_{1 \le m \le s} \hat{w}_m^{1/2} \sum_{i=1}^s \sum_{j=1}^s \tau |(U_j, \hat{w}_i^{1/2}V_i)|
      \\
      &\le \max_{1 \le \ell,k \le s}[(\mathbb{A}_r^{-1})_{\ell k}] \frac{\max_{1 \le m \le s} \hat{w}_m^{1/2}}{\min_{1 \le n \le s} \hat{w}_n^{1/2}} \sum_{i=1}^s \norm[0]{(\tau \hat{w}_i)^{1/2}U_i}_{L^2(\Omega)} \sum_{j=1}^s \norm[0]{(\tau \hat{w}_j)^{1/2} V_j}_{L^2(\Omega)}
      \\
      &= \max_{1 \le \ell,k \le s}[(\mathbb{A}_r^{-1})_{\ell k}] \frac{\max_{1 \le m \le s} \hat{w}_m^{1/2}}{\min_{1 \le n \le s} \hat{w}_n^{1/2}} \norm[0]{\boldsymbol{u}_h}_{r,\boldsymbol{v}} \norm[0]{\boldsymbol{v}_h}_{r,\boldsymbol{v}}.
    \end{split}
  \end{equation*}
  Next, again using Cauchy--Schwarz,
  \begin{equation*}
    T_2
    \le \sum_{i=1}^s \norm[0]{(\hat{w}_i\nu)^{1/2} \nabla U_i}_{L^2(\Omega)} \norm[0]{(\hat{w}_i\nu)^{1/2} \nabla V_i}_{L^2(\Omega)}
    =    
    \norm[0]{\boldsymbol{u}_h}_{r,\boldsymbol{v}} \norm[0]{\boldsymbol{v}_h}_{r,\boldsymbol{v}}.
  \end{equation*}
  The bounds for $T_1$ and $T_2$, together with
  \cref{eq:aboundp-wradauiia}, imply \cref{eq:boldabound-wradauiia}. To
  prove \cref{eq:boldbbound-wradauiia}, note that
  \begin{multline*}
    |\boldsymbol{b}_r(\boldsymbol{v}_h, \boldsymbol{q}_h)|
    = \envert[3]{\sum_{i=1}^s \hat{w}_i b(V_i, Q_i)}
    \le \sum_{i=1}^s \hat{w}_i^{1/2} \norm[0]{V_i}_v \hat{w}_i^{1/2} \norm[0]{Q_i}_q
    \\
    \le \del[3]{ \sum_{i=1}^s \hat{w}_i \norm[0]{V_i}_v^2 }^{1/2}
    \del[3]{ \sum_{i=1}^s \hat{w}_i \norm[0]{Q_i}_q^2 }^{1/2}
    = \norm[0]{\boldsymbol{v}_h}_{r,\boldsymbol{v}} \norm[0]{\boldsymbol{q}_h}_{r,\boldsymbol{q}},      
  \end{multline*}
  where the first inequality was proven in \cite[Section
  2]{mardal2004uniform} and the second by using the discrete
  H\"older's inequality.
\end{proof}

\subsection{Discontinuous Galerkin time-stepping}

We now consider DG time-stepping to discretize \cref{eq:stokes} in
time. In \cref{sss:wspacetime} we first analyze this space-time
DG+Taylor--Hood discretization. Then, in \cref{sss:wspacetimeRK}, we
derive a new RK time-stepping scheme that is algebraically equivalent
to the space-time scheme.

\subsubsection{Space-time discretization}
\label{sss:wspacetime}

On each time interval we define the space-time finite element spaces
$\mathcal{V}_h^n := \mathcal{P}_{k_t}(I_n;V_h)$ and
$\mathcal{Q}_h^n := \mathcal{P}_{k_t}(I_n;Q_h)$, where
$\mathcal{P}_k(I_n;X)$ denotes the set of polynomials of degree $k_t$
in time with values in $X$ (see, for example,
\cite{chrysafinos2010discontinuous}). We use weighted inner products
in time to define the fully discrete time-dependent Stokes
problem. Weighted inner products have previously been used to define
space-time discretizations in, for example,
\cite{french1993space,sosa2020space}.

Let $\kappa_n(t) := \exp(-\tau(t-t_n))$. The discretization is given
by: For $n=0,\hdots,N-1$, find $u_h^n \in \mathcal{V}_h^n$ and
$p_h^n\in \mathcal{Q}_h^n$ such that
\begin{equation}
  \label{eq:discreteStokes}
  \boldsymbol{\hat{a}}(u_h^n,v_h) + \boldsymbol{\hat{b}}(v_h,p_h^n) - \boldsymbol{\hat{b}}(u_h^n,q_h)
  = \int_{I_n} (f, v_h) \kappa_n \dif t + (u_n^-, v_n^+), 
\end{equation}
for all $v_h \in \mathcal{V}_h^n$ and $q_h \in \mathcal{Q}_h^n$, where
\begin{equation*}
  \begin{split}
    \boldsymbol{\hat{a}}(u_h^n,v_h)
    &:= \int_{I_n} (\partial_t u_h^n, v_h) \kappa_n \dif t
      + (u_n^+, v_n^+)
      + \int_{I_n} a(u_h^n, v_h) \kappa_n \dif t,
    \\
    \boldsymbol{\hat{b}}(v_h,q_h)
    &:= \int_{I_n} b(v_h, q_h) \kappa_n \dif t,    
  \end{split}
\end{equation*}
in which $u_n^- := u_h^{n-1}(t_n)$ and $u_n^+ := u_h^n(t_n)$. Note
that if $\kappa_n(t) \equiv 1$ then we recover the discretization
analyzed in \cite{chrysafinos2010discontinuous}.

To show that $\boldsymbol{\hat{a}}(\cdot, \cdot)$ and
$\boldsymbol{\hat{b}}(\cdot, \cdot)$ satisfy
\cref{eq:aconditions,eq:bconditions}. We introduce the following inner
products on $\mathcal{V}_h^n$ and $\mathcal{Q}_h^n$:
\begin{align*}
  (u_h,v_h)_{v,t} & := \int_{I_n} (u_h,v_h)_v \kappa_n \dif t,
  &
  (p_h,q_h)_{q,t} & := \int_{I_n}(p_h,q_h)_q \kappa_n \dif t,
\end{align*}
and denote their induced norms by $\norm[0]{\cdot}_{v,t}$ and
$\norm[0]{\cdot}_{q,t}$, respectively. Here $(\cdot, \cdot)_v$ and
$(\cdot, \cdot)_q$ are defined in \cref{eq:innerprodvq}.

\begin{lemma}[Coercivity]
  \label{lem:st-coercivity}
  For all $v_h \in \mathcal{V}_h^n$,
  \begin{equation*}
    \boldsymbol{\hat{a}}(v_h, v_h)
    \ge
    \tfrac{1}{2} \norm[0]{v_h}_{v,t}^2.
  \end{equation*}
\end{lemma}
\begin{proof}
  We have, using integration by parts and that
  $\kappa_n' = -\tau \kappa_n$,
  \begin{align*}
    \boldsymbol{\hat{a}}(v_h, v_h)
    =& \int_{I_n} \tfrac{1}{2} \od{}{t} \norm[0]{v_h}^2 \kappa_n \dif t
       + \norm[0]{v_n^+}^2 + \int_{I_n} a(v_h, v_h) \kappa_n \dif t
    \\
    \ge & \tfrac{1}{2} \int_{I_n} \tau\norm[0]{v_h}^2 \kappa_n \dif t
      + \int_{I_n} \nu \norm[0]{\nabla v_h}^2 \kappa_n \dif t, 
  \end{align*}
  which is the desired result.
\end{proof}

\begin{lemma}[Stability]
  \label{lem:st-stability}
    There exists a uniform constant $\beta > 0$ such that
    \begin{equation*}
      \sup_{v_h \in \mathcal{V}_h} \frac{\boldsymbol{\hat{b}}(v_h,q_h)}{\norm[0]{v_h}_{v,t}}
      \ge
      \beta \norm[0]{q_h}_{q,t}
      \qquad
      \forall q_h \in \mathcal{Q}_h.
    \end{equation*}
\end{lemma}
\begin{proof}
  Following the proof of \cite[Theorem 3.2]{kirk2023analysis}, we use
  the tensor-product structure of $\mathcal{Q}_h$ to write
  $q_h \in \mathcal{Q}_h$ in terms of an orthonormal basis
  $\cbr[0]{\phi_i}_{i=0}^{k_t}$ of $P_k(I_n)$ with respect to the
  weighted $L^2(I_n)$ inner product,
  \begin{equation*}
      q_h(x,t) = \sum_{i=0}^{k_t} \phi_i(t) q_i(x), \quad \phi_i \in P_{k_t}(I_n),\ q_i \in Q_h.
  \end{equation*}
  By \cref{eq:infsupmardal}, for a given $q_i \in Q_h$ there exists
  $\tilde{v}_i \in V_h$ such that
  $b(\tilde{v}_i, q_i) = \norm[0]{q_i}_q^2$ and
  $\norm[0]{\tilde{v}_i}_v \le \beta^{-1} \norm[0]{q_i}_q$. We now
  define
  $\tilde{v}_h(x,t) = \sum_{i=0}^{k_t} \phi_i(t) \tilde{v}_i \in
  \mathcal{V}_h$. Given $q_h \in \mathcal{Q}_h$, and that
  $\cbr[0]{\phi_i}$ is an orthonormal basis with respect to the
  weighted $L^2(I_n)$ inner product, we find
  \begin{equation}
    \label{eq:infsupeq}
    \norm[0]{q_h}_{q,t}^2
    = \sum_{i=0}^{k_t} \norm[0]{q_i}_q^2
    = \sum_{i=0}^{k_t} b(\tilde{v}_i,q_i)
    = \boldsymbol{\hat{b}}(\tilde{v}_h,q_h).
  \end{equation}
  Likewise, we find that
  \begin{equation}
    \label{eq:infsupineq}
    \norm[0]{\tilde{v}_h}_{v,t}^2
    = \sum_{i=0}^{k_t}\norm[0]{\tilde{v}_i}_v^2
    \le \beta^{-2} \sum_{i=0}^{k_t}\norm[0]{\tilde{q}_i}_q^2
    = \beta^{-2} \norm[0]{q_h}_{q,t}^2.
  \end{equation}
  The result follows from \cref{eq:infsupeq,eq:infsupineq}.
\end{proof}

\begin{lemma}[Continuity]
  \label{lem:st-continuity}
  There exist uniform constants $c_a, c_b > 0$ such that
  \begin{subequations}
    \begin{align}
      \label{eq:st-continuity-a}
      \envert[0]{\boldsymbol{\hat{a}}(u_h,v_h)}
      & \le c_a \norm[0]{u_h}_{v,t} \norm[0]{v_h}_{v,t} && \forall u_h, v_h \in \mathcal{V}_h^n,
      \\
      \label{eq:st-continuity-b}
      \envert[0]{\boldsymbol{\hat{b}}(v_h, q_h)}
      & \le c_b \norm[0]{v_h}_{v,t} \norm[0]{q_h}_{v,t} && \forall v_h \in \mathcal{V}_h^n,\ q_h \in \mathcal{Q}_h^n.
    \end{align}
  \end{subequations}
\end{lemma}
\begin{proof}
  Using the triangle inequality,
  \begin{equation*}
    \envert[0]{\boldsymbol{\hat{a}}(u_h,v_h)}
    \leq
    \envert[3]{\int_{I_n} (\partial_t u_h, \kappa_n v_h) \dif t }
    + \envert[3]{ (u_n^+, v_n^+) }
    + \envert[3]{ \int_{I_n} a(u_h, v_h) \kappa_n \dif t }
    =: T_1 + T_2 + T_3.
  \end{equation*}
  We bound $T_1$, $T_2$, and $T_3$ separately. Note that
  $e^{-1} \le \kappa_n(t) \le 1$ for $t \in I_n$. Using the
  Cauchy--Schwarz inequality and that
  $\int_{I_n} \norm[0]{\partial_tv_h}^2 \dif t \le c_t \tau \int_{I_n}
  \norm[0]{v_h}^2 \dif t$ for all $v_h \in \mathcal{V}_h$, where
  $c_t > 0$ is a uniform constant (see, e.g.,
  \cite[eq. (19b)]{kirk2023analysis}), we find
  \begin{equation}
    \label{eq:st-continuity-aux1}
    \begin{split}
      T_1^2
      & \le \int_{I_n} \kappa_n \norm[0]{\partial_t u_h}^2 \dif{t} \int_{I_n} \kappa_n \norm[0]{ v_h}^2 \dif{t}
      \\
      & \le c_t^2 e \int_{I_n} \tau  \kappa_n \norm[0]{u_h}^2 \dif{t} \int_{I_n} \tau  \kappa_n \norm[0]{ v_h}^2 \dif{t}
      \le c_t^2 e \norm[0]{u_h}_{v,t}^2 \norm[0]{v_h}_{v,t}^2.
    \end{split}
  \end{equation}
  Next, again using the Cauchy--Schwarz inequality but now in
  combination with
  $\norm[0]{u_n^+}^2 \le c_{\tau} \tau \int_{I_n} \norm[0]{u_h}^2 \dif
  t$, which holds for all $u_h \in \mathcal{V}_h$ and where
  $c_{\tau}>0$ is a uniform constant (see, e.g., \cite[Lemma 4,
  eq. (17b)]{sosa2020space}), 
  \begin{equation}
    \label{eq:st-continuity-aux2}
      T_2^2
      \le \norm[0]{u_h^+}^2 \norm[0]{v_h^+}^2
      \le c_{\tau}^2 e^2 \int_{I_n} \tau \kappa_n \norm[0]{u_h}^2 \dif{t} \int_{I_n} \tau \kappa_n \norm[0]{v_h}^2 \dif{t}
      \le c_{\tau}^2 e^2 \norm[0]{u_h}_{v,t}^2 \norm[0]{v_h}_{v,t}^2.
  \end{equation}
  To bound $T_3$ we use the Cauchy--Schwarz inequality and that
  $|a(w_h,z_h)| \le \nu \norm[0]{\nabla w_h}_v\norm[0]{\nabla z_h}_v$ for
  all $w_h,z_h \in V_h$ to find
  \begin{equation}
    \label{eq:st-continuity-aux3}
    T_3 
    \le \int_{I_n} \kappa_n^{1/2} \norm[0]{u_h}_v \kappa_n^{1/2} \norm[0]{v_h}_v \dif{t} 
    \le \norm[0]{u_h}_{v,t} \norm[0]{v_h}_{v,t}.
  \end{equation}
  \Cref{eq:st-continuity-a} follows after combining
  \cref{eq:st-continuity-aux1,eq:st-continuity-aux2,eq:st-continuity-aux3}. Next,
  it is shown in \cite[Section 2]{mardal2004uniform} that
  $|b(w_h,r_h)| \le \norm[0]{w_h}_v \norm[0]{r_h}_q$ for
  all $w_h \in V_h$ and $r_h \in Q_h$. Using also Cauchy--Schwarz, we
  find
  \begin{equation*}
    \envert[0]{\boldsymbol{\hat{b}}(v_h, q_h)}
    \le \int_{I_n} \kappa_n^{1/2} \norm[0]{v_h}_v \kappa_n^{1/2} \norm[0]{q_h}_q \dif{t}
    \le \norm[0]{v_h}_{v,t} \norm[0]{q_h}_{v,t},
  \end{equation*}
  proving \cref{eq:st-continuity-b}.
\end{proof}

\subsubsection{A new RK time-stepping method}
\label{sss:wspacetimeRK}

Traditional DG time-stepping methods (for which $\kappa_n=1$) can be
written as a collocation method resulting in a Radau IIA RK method,
see for example, \cite{Akrivis:2022,Akrivis:2011}. The bilinear form
$\boldsymbol{\widetilde{a}}_r$ of the Radau IIA RK method
\cref{eq:stokes-radauiiaTH}, however, cannot be shown to satisfy
\cref{eq:abs-coercivity}. If there exists a matrix $\widehat{W}$ that
satisfies the conditions of \cref{lem:coercivity-wradauiia}, then the
bilinear form $\boldsymbol{a}_r$ of the stage-weighted Radau IIA RK
method \cref{eq:stokes-wradauiia} is coercive.

An alternative to using a stage-weighted Radau IIA time discretization
\cref{eq:stokes-wradauiia} is to use an RK time discretization such
that the resulting discretization is algebraically equivalent to
\cref{eq:discreteStokes}. This because, as shown in
\cref{sss:wspacetime}, the space-time discretization
\cref{eq:discreteStokes} results in bilinear forms that satisfy
\cref{eq:aconditions,eq:bconditions}. In this section we derive this
new RK method.

Let $\hat{t} = \tau (t-t_n) \in [0,1]$ and
$\hat{\kappa}_n(\hat{t}) = \exp(-\hat{t})$. Since
$(1-\hat{t}) \hat{\kappa}_n(\hat{t})$ is positive on $(0,1)$ there
exists a unique monic polynomial of degree $(s-1)$,
$\phi_{s-1}(\hat{t})$, such that
\begin{equation}
  \label{eq:orthogonalityphismin1}
  \int_0^1 \phi_{s-1}(\hat{t}) q(\hat{t}) (1-\hat{t}) \hat{\kappa}_n(\hat{t}) \dif \hat{t} = 0
  \qquad \forall q \in \mathcal{P}_{s-2}([0,1]).
\end{equation}
For the roots of $\phi_{s-1}(\hat{t})$,
$\cbr[0]{c_{i,\tau}}_{i=1}^{s-1}$, let
\begin{equation}
  \label{eq:weighted-quadrature}
  0 < c_{1,\tau} < \hdots < c_{s-1,\tau}< c_{s,\tau} = 1.
\end{equation}
For Lagrange basis and weights 
\begin{equation}
  \label{eq:lagrange-and-weight}
  \ell_i(\hat t)=\prod_{\substack{j=1\\j\ne i}}^s \frac{\hat t-c_{j,\tau}}{c_{i,\tau}-c_{j,\tau}},
  \qquad
  w_i=\int_0^1\ell_i(\hat t)\hat\kappa_n(\hat t)\dif\hat t,
  \quad i=1,\ldots,s,
\end{equation}
we define a weighted Radau quadrature rule for $q \in C([0,1])$ by 
\begin{equation}
  \label{eq:weighted-quadrature-rule}
  \sum_{i=1}^s w_i q(c_{i,\tau}) .
\end{equation}
We refer to \cite{trostorff2016} for the existence of weighted Radau
quadrature.

\begin{lemma}
  For any $p \in \mathcal{P}_{2s-2}([0,1])$ the quadrature rule
  \cref{eq:weighted-quadrature-rule} satisfies
  \begin{equation}
    \label{eq:weighted-radau-exactness}
    \int_0^1p(\hat t)\hat\kappa_n(\hat t)\dif\hat t
    =\sum_{i=1}^s w_i p(c_{i,\tau})
    \qquad \forall p\in\mathcal{P}_{2s-2}([0,1]).
  \end{equation}
  Furthermore, $w_i >0$ for $i=1, \ldots, s$.
\end{lemma}
\begin{proof}
  Any $v\in\mathcal{P}_{s-1}([0,1])$ can be written in terms of the
  Lagrange basis:
  \begin{equation}
    \label{eq:lagrange-interpolation}
    v(\hat t)=\sum_{i=1}^s v(c_{i,\tau})\ell_i(\hat t).
  \end{equation}
  Multiplying \cref{eq:lagrange-interpolation} by $\hat\kappa_n$ and
  integrating proves that the quadrature formula and the integral
  coincide for $v$.

  There exist polynomials $q\in\mathcal{P}_{s-2}([0,1])$ and
  $v\in\mathcal{P}_{s-1}([0,1])$ such that a given
  $p\in\mathcal{P}_{2s-2}([0,1])$ can be written as
  \begin{equation}
    \label{eq:polynomial-division}
    p(\hat t)=(\hat t-1)\phi_{s-1}(\hat t)q(\hat t)+v(\hat t).
  \end{equation}
  Using orthogonality of $\phi_{s-1}$ (see
  \cref{eq:orthogonalityphismin1}), applying the exact integration of
  the quadrature formula when applied to $v$, and that the values of
  $v$ and $p$ agree at every quadrature node (because the first term
  of \cref{eq:polynomial-division} vanishes at the roots of
  $\phi_{s-1}$ and at $1$), gives
  \begin{equation*}
    \int_0^1p(\hat t)\hat\kappa_n (\hat t) \dif\hat t
    =\int_0^1v (\hat t)\hat\kappa_n (\hat t)\dif\hat t
    =\sum_{i=1}^s w_i v(c_{i,\tau})
    =\sum_{i=1}^s w_i p(c_{i,\tau}),    
  \end{equation*}
  which proves \cref{eq:weighted-radau-exactness}. That $w_i>0$,
  $i=1,\hdots s$, follows from this exact integration because
  $\int_0^1 \ell_i(\hat t)^2 \hat{\kappa}_n(\hat{t}) \dif \hat{t} =
  \sum_{j=1}^s w_j \ell_i(c_{j,\tau})^2 = w_i$.
\end{proof}

\begin{lemma}
  \label{lem:HHinvW}
  For $1\le i,j\le s$, $s \ge 2$, let
  \begin{equation}
    \label{eq:temporal-matrices}
    D_{ij}:=\ell_j'(c_{i,\tau}),\quad d_i:=\ell_i(0),\quad  W:=\text{diag}(w_1,\ldots,w_s),\quad H:=WD+dd^T,
  \end{equation}
  where the prime in the definition of $D_{ij}$ denotes
  differentiation with respect to $\hat t$.  Then, $H$ is invertible
  and the following hold:
  \begin{equation}
    \label{eq:tableau-identities}
    H\mathbf1=d,\qquad H^{-1}W \mathbf1=c_\tau,  \qquad e_s^T H^{-1}W \mathbf1=1,
  \end{equation}
  where $\mathbf1=(1,\ldots,1)^T\in\mathbb{R}^s$,
  $e_s = (0, \ldots, 0,1)^T$,
  $c_{\tau}=(c_{1,\tau}, \ldots, c_{s,\tau})^T$.
\end{lemma}
\begin{proof}
  We first show that $H$ is invertible by showing that if
  $z \in \mathbb{R}^s$ is such that $Hz=0$, then $z=0$. For
  $z \in \mathbb{R}^s$ set
  $p_z(\hat t) = \sum_{j=1}^s z_j \ell_j(\hat t) \in
  \mathcal{P}_{s-1}([0,1])$. Then,
  \begin{align*}
    &\int_0^1 \hat{\kappa}_n p_z' (\hat{t}) \ell_i(\hat{t}) \dif \hat{t} + p_z(0) \ell_i(0)
    \\
    &= \sum_{k=1}^s w_k \left[ \sum_{j=1}^s z_j \ell_j'(c_{k,\tau}) \ell_i(c_{k,\tau}) \right]+ \ell_i(0) \sum_{j=1}^s z_j \ell_j(0) \quad (\text{by~\cref{eq:weighted-radau-exactness}}) 
    \\
    &= \sum_{k=1}^s w_k \left[\ell_i(c_{k,\tau}) \sum_{j=1}^s D_{kj} z_j \right] + d_i d^T z  \quad (\text{by~\cref{eq:temporal-matrices}})
    \\
    &= w_i (Dz)_i + d_i d^T z \quad (\text{by~\cref{eq:temporal-matrices} and }\ell_i(c_{k,\tau}) = \delta_{ik})
    \\
    &= (Hz)_i. \quad (\text{by~\cref{eq:temporal-matrices}})
  \end{align*}
  Suppose now that $Hz = 0$ for some $z \in \mathbb{R}^s$. Since any
  $v \in \mathcal{P}_{s-1}([0,1])$ can be written as a linear
  combination of $\ell_i$'s, the above identity gives
  \begin{equation*}
    \int_0^1 \hat{\kappa}_n (\hat t) p_z'(\hat t) v(\hat t) \dif \hat{t} + p_z(0) v(0) = 0
    \qquad \forall v \in \mathcal{P}_{s-1}([0,1]).
  \end{equation*}
  Choosing
  $v(\hat t) = \hat{t} p_z'(\hat t) \in \mathcal{P}_{s-1}([0,1])$
  implies
  $\int_0^1 \hat{t} \hat{\kappa}_n(\hat t)|p_z'(\hat{t})|^2 \dif
  \hat{t} = 0$, and so $p_z' \equiv 0$. Moreover, choosing
  $v(\hat{t})\equiv 1$ implies $p_z(0) = 0$, and so $p_z \equiv
  0$. Therefore, $z = 0 \in \mathbb{R}^s$ because
  $z_i = p_z(c_{i,\tau})$.

  We now show the identities in \cref{eq:tableau-identities}. Note
  that Lagrange interpolation of the constant function $1$ and
  $\hat t$ gives
  \begin{equation*}
    \sum_{j=1}^s\ell_j(\hat t)=1, \qquad \sum_{j=1}^s c_{j,\tau}\ell_j(\hat t)=\hat t.
  \end{equation*}
  Differentiating these identities and evaluating at the nodes, and
  also evaluating them at $\hat{t}=0$, yields
  \begin{equation*}
    D\mathbf1=0,
    \quad
    d^T\mathbf1=1,
    \qquad
    Dc_\tau=\mathbf1,
    \quad d^Tc_\tau=0.    
  \end{equation*}
  From these identities and the definition of $H$,
  \begin{equation*}
    H\mathbf1=d,\qquad Hc_\tau=W\mathbf1,
  \end{equation*}
  from which the first and second identity in
  \cref{eq:tableau-identities} follow. The third identity follows from
  $e_s^T H^{-1}W \mathbf1=e_s^T c_{\tau}=1$, since $c_{s,\tau}=1$.
\end{proof}

\begin{lemma}
  \label{lem:weighted-RK}
  Let $s=k_t+1 \ge 2$, $t_{n,k} = t_n + c_{k,\tau} \Delta t$ with
  quadrature nodes
  $0 < c_{1,\tau} < \hdots < c_{s-1,\tau}< c_{s,\tau} = 1$ as given by
  \cref{eq:weighted-quadrature}. Let
  $\mathbb{A} = [a_{ij}] = H^{-1}W$. The space-time discretization
  \cref{eq:discreteStokes} of the time-dependent Stokes equations is
  equivalent to using RK time-stepping with Butcher tableau
  \begin{equation}
    \label{eq:butcher-tableau}
    \begin{array}{c|c}
      c_\tau&\mathbb A\\ \hline
            & e_s^T\mathbb A
    \end{array},
    \qquad e_s=(0,\ldots,0,1)^T\in\mathbb{R}^s,
  \end{equation}
  if 
  \begin{equation}
    \label{eq:load-assumption}
    \int_{I_n} \kappa_n (t) (f(t), v_h^n(t)) \dif t = \Delta t \sum_{i=1}^s w_i (f(t_{n,i}), v_h^n(t_{n,i}))
    \quad \forall v_h^n \in \mathcal{V}_h^n.  
  \end{equation}
\end{lemma}
\begin{proof}[Proof of \cref{lem:weighted-RK}]
  Let $t_{n,i}:=t_n+c_{i,\tau}\Delta t$ and write the Lagrange
  expansion of $u_h^n$, $p_h^n$ on the physical slab as
  \begin{equation}
    \label{eq:time-nodal-expansions}
    u_h^n(t)=\sum_{j=1}^sU_j\ell_j(\hat t),
    \qquad
    p_h^n(t)=\sum_{j=1}^sP_j\ell_j(\hat t),
    \qquad \hat t=\frac{t-t_n}{\Delta t}
  \end{equation}
  where $U_j\in V_h$, $P_j\in Q_h$. Suppose that $(U_j, P_j)_{j=1}^s$
  are obtained by the DG method. It follows that
  \begin{equation}
    \label{eq:derivative-and-left-trace}
    \partial_tu_h^n(t_{n,i})
    =\frac{1}{\Delta t}\sum_{j=1}^sD_{ij}U_j,
    \qquad u_n^+=\sum_{j=1}^sd_jU_j.
  \end{equation}
  In \cref{eq:discreteStokes}, fix $v_h\in V_h$, $q_h \in Q_h$, and
  use the time-dependent test function $v_h\ell_i(\hat t)$,
  $q_h\ell_i(\hat{t})$. The time derivative product has temporal
  degree at most $(s-2)+(s-1)=2s-3$, and the velocity and pressure
  products have temporal degree at most $2s-2$. Thus all integrals on
  $I_n$ are evaluated exactly by the weighted quadrature given in
  \cref{eq:weighted-quadrature-rule}. More specifically, applying the
  weighted quadrature rule to the derivative term gives
  \begin{align*}
    \int_{I_n}\kappa_n(t)(\partial_tu_h^n(t),v_h) \ell_i(\hat t)\dif t
    &= \Delta t\sum_{k=1}^s w_k (\partial_tu_h^n(t_{n,k}),v_h)\ell_i(c_{k,\tau})
    \\
    &= w_i\sum_{j=1}^sD_{ij}(U_j,v_h)
      \quad (\text{by \cref{eq:derivative-and-left-trace} and }\ell_i(c_{k,\tau}) = \delta_{ik}).
  \end{align*}
  The term corresponding to $(u_n^+, v_n^+)$ can be written as
  $d_i \sum_{j=1}^s d_j (U_j,v_h)$. Applying the weighted quadrature
  to \cref{eq:discreteStokes} therefore results in
  \begin{multline}
    \label{eq:nodal-momentum}
    \sum_{j=1}^sH_{ij}(U_j,v_h)
    +\Delta t w_i\bigl[a(U_i,v_h)+b(v_h,P_i)+b(U_i,q_h)\bigr]
    \\
    =d_i (u_n^-,v_h)+\Delta t w_i (f(t_{n,i}),v_h), \qquad i=1,\hdots,s.
  \end{multline}
  If $v_h=0$, then $b(U_i, q_h)=0$ for all $q_h \in Q_h$ because
  $\Delta t w_i b(U_i, q_h) = 0$ and $w_i>0$.  If we multiply
  \cref{eq:nodal-momentum} for $i=1,\ldots,s$ by $H^{-1}$ in the stage
  index, use that $H^{-1}d=\mathbf1$, $H^{-1}W=\mathbb A$, and define
  $K_j \in V_h$ by
  \begin{equation*}
    (K_j, v_h) = (f(t_{n,j}), v_h) - a(U_j, v_h) - b(v_h, P_j) \qquad \forall v_h \in V_h,
  \end{equation*}
  we find that $U_i = u_n^- + \Delta t \sum_{j=1}^s a_{ij} K_j$. In
  particular,
  $U_s = u_h^n(t_{n+1}^-) = u_n^- + \Delta t \sum_{j=1}^s a_{sj}K_j$.
  By recalling that the last row of the Butcher tableau is
  $e_s^{T} \mathbb{A}$, $u_h^n(t_{n+1}^-)$ is the solution of the RK
  method. Conversely, let \((U_i,P_i)_{i=1}^s\) satisfy the RK stage
  equations and reconstruct \(u_h^n\) and \(p_h^n\) by
  \(\cref{eq:time-nodal-expansions}\). Multiplying the stage momentum
  equations by \(H\), and using \(H\mathbf1=d\) and \(H\mathbb A=W\),
  recovers the momentum part of \(\cref{eq:nodal-momentum}\). The
  stage constraints give its remaining part. By exactness of the
  weighted quadrature and the assumption \cref{eq:load-assumption},
  these are the DG equations tested with \(v_h\ell_i\) and
  \(q_h\ell_i\). Since these test functions span the DG test spaces,
  the reconstructed pair solves the DG problem.
\end{proof}

As proven in \cref{lem:weighted-RK} we can write
\cref{eq:discreteStokes} as \cref{eq:nodal-momentum}. Using the
notation $f_i(x) = f(x,t_{n,i})$, that
$w_i^{-1}H_{ij} = (W^{-1}H)_{ij} = (\mathbb{A}^{-1})_{ij}$ and, thanks
to \cref{eq:tableau-identities},
$w_i^{-1}d_i = (W^{-1}d)_i = (W^{-1}H\mathbf{1})_i =
(\mathbb{A}^{-1}\mathbf{1})_i$, we can write \cref{eq:nodal-momentum}
as:
\begin{multline*}
  w_i \tau \sum_{j=1}^s(\mathbb{A}^{-1})_{ij}(U_j,v_h)
  + w_i \bigl[a(U_i,v_h)+b(v_h,P_i)+b(U_i,q_h)\bigr]
  \\
  = \tau w_i (\mathbb{A}^{-1}\mathbf{1})_i (u_n^-,v_h) + w_i (f_i,v_h),
  \qquad i=1,\hdots,s.
\end{multline*}
Advancing from $t_n$ to $t_{n+1}$ we consider the following fully
discrete problem: find
$\boldsymbol{u}_h := (U_1,\hdots,U_s) \in \boldsymbol{V}_h$ and
$\boldsymbol{p}_h := (P_1,\hdots,P_s) \in \boldsymbol{Q}_h$ such that
\begin{equation}
  \label{eq:stokesradauw}
  \boldsymbol{a}(\boldsymbol{u}_h, \boldsymbol{v}_h)
  + \boldsymbol{b}(\boldsymbol{v}_h, \boldsymbol{p}_h)
  - \boldsymbol{b}(\boldsymbol{u}_h, \boldsymbol{q}_h)
  = \sum_{i=1}^s w_i \sbr[3]{(f_i, V_i)
  + \tau \sum_{j=1}^s (\mathbb{A}^{-1})_{ij}(u_n^{-}, V_i)},
\end{equation}
for all $\boldsymbol{v}_h := (V_1,\hdots,V_s) \in \boldsymbol{V}_h$
and $\boldsymbol{q}_h := (Q_1,\hdots,Q_s) \in \boldsymbol{Q}_h$, in
which
\begin{equation*}
  \boldsymbol{a}(\boldsymbol{u}_h, \boldsymbol{v}_h)
  := \sum_{i=1}^s w_i \sbr[3]{\tau \sum_{j=1}^s (\mathbb{A}^{-1})_{ij}(U_j, V_i)
    +  a(U_i, V_i)},
  \quad
  \boldsymbol{b}(\boldsymbol{v}_h, \boldsymbol{p}_h)
  := \sum_{i=1}^s w_i b(V_i, P_i).
\end{equation*}

We will show that $\boldsymbol{a}(\cdot, \cdot)$ and
$\boldsymbol{b}(\cdot, \cdot)$ satisfy
\cref{eq:aconditions,eq:bconditions}. We introduce the following inner
products associated with $\boldsymbol{V}_h$ and $\boldsymbol{Q}_h$:
\begin{equation*}
  (\boldsymbol{u}_h, \boldsymbol{v}_h)_{\boldsymbol{v}} := \sum_{i=1}^s
  w_i (U_i, V_i)_v,
  \qquad
  (\boldsymbol{p}_h, \boldsymbol{q}_h)_{\boldsymbol{q}} := \sum_{i=1}^s
  w_i (P_i, Q_i)_q,
\end{equation*}
respectively, with induced norms denoted by
$\norm[0]{\cdot}_{\boldsymbol{v}}$ and
$\norm[0]{\cdot}_{\boldsymbol{q}}$. Here $(\cdot, \cdot)_v$ and
$(\cdot, \cdot)_q$ are defined in \cref{eq:innerprodvq}.

Proving stability of $\boldsymbol{b}$ and continuity of
$\boldsymbol{a}$ and $\boldsymbol{b}$ follows identical steps as used
to prove \cref{lem:stability-wradauiia,lem:continuity-wradauiia}. The
following lemma will be used to show coercivity of
$\boldsymbol{a}(\cdot, \cdot)$.

\begin{lemma}
  \label{lem:WAinvWAinvTspd}
  Suppose that $\mathbb{A}$ is the Butcher tableau defined in
  \cref{lem:weighted-RK}. Let $W$ be defined in
  \cref{eq:temporal-matrices}. Then
  $\tfrac{1}{2}( W\mathbb{A}^{-1} + (W\mathbb{A}^{-1})^T)$ is
  symmetric positive definite.
\end{lemma}
\begin{proof}
  Recall that $\mathbb{A} = H^{-1}W$, so
  $\tfrac{1}{2}( W\mathbb{A}^{-1} + (W\mathbb{A}^{-1})^T) =
  \tfrac{1}{2}(H + H^{T})$. We show that $H + H^{T}$ is symmetric
  positive definite. From \cref{eq:temporal-matrices} and the exact
  integration of the quadrature rule,
  \begin{align*}
    (WD)_{ij}
    &= w_i \ell_j'(c_{i,\tau}) = \sum_{k=1}^s w_k \ell_i(c_{k,\tau}) \ell_j'(c_{k,\tau})
      \quad (\because \ell_i(c_{k,\tau}) = \delta_{ik})
    \\
    &= \int_0^1 \hat{\kappa}_n(\hat{t}) \ell_i(\hat{t}) \ell_j'(\hat{t}) \dif \hat{t} .
  \end{align*}
  Then,
  \begin{equation*}
    (WD+D^TW)_{ij}
    = \int_0^1 \hat{\kappa}_n(\hat{t}) [\ell_i(\hat{t}) \ell_j'(\hat{t}) + \ell_j(\hat{t}) \ell_i'(\hat{t}) ] \dif \hat{t}
    = \int_0^1 \hat{\kappa}_n(\hat{t}) (\ell_i(\hat{t}) \ell_j(\hat{t}))'\dif \hat{t}.
  \end{equation*}
  By the integration by parts and
  $\hat{\kappa}_n'(\hat{t}) = -\hat{\kappa}_n(\hat{t})$,
  \begin{align*}
    (WD+D^TW)_{ij}
    &= (\hat{\kappa}_n \ell_i \ell_j) |_0^1 + \int_0^1 \hat{\kappa}_n (\hat{t}) \ell_i(\hat{t}) \ell_j(\hat{t}) \dif \hat{t}
    \\
    &= e^{-1} \ell_i(1) \ell_j(1) - \ell_i(0) \ell_j(0) + \int_0^1 \hat{\kappa}_n (\hat{t}) \ell_i(\hat{t}) \ell_j(\hat{t}) \dif \hat{t}.
  \end{align*}
  Applying the quadrature rule to the last integral, a matrix form of
  this identity is
  \begin{equation*}
    WD + D^T W = e^{-1} e_s e_s^T - d d^T + W .
  \end{equation*}
  The result follows by
  \begin{equation*}
    H + H^T
    = (WD + d d^T) + (D^T W + dd^T)
    = d d^T + e^{-1}e_s e_s^T + W.
  \end{equation*}
\end{proof}
\begin{lemma}[Coercivity]
  \label{lem:coercivity-newradau}
  There exists a uniform constant $\alpha > 0$ such that
  \begin{equation*}
    \boldsymbol{a}(\boldsymbol{v}_h, \boldsymbol{v}_h) \ge \alpha \norm[0]{\boldsymbol{v}_h}_{\boldsymbol{v}}^2
    \qquad \forall \boldsymbol{v}_h \in \boldsymbol{V}_h.
  \end{equation*}
\end{lemma}
\begin{proof}
  Since $C = \tfrac{1}{2}( W\mathbb{A}^{-1} + (W\mathbb{A}^{-1})^T)$
  is symmetric positive definite by \cref{lem:WAinvWAinvTspd}, the
  proof follows the same steps as that of
  \cref{lem:coercivity-wradauiia}.
\end{proof}

\section{Preconditioning}
\label{s:preconditioning}

Let $\mathrm{u},\mathrm{v}$ be the expansion coefficients related to
$u_h,v_h \in V_h$ and let $\mathrm{p},\mathrm{q}$ be the expansion
coefficients related to $p_h,q_h \in Q_h$, respectively. Let $P_v$,
$P_q$, and $B$ be the matrices defined by
$\langle P_v \mathrm{u}, \mathrm{v} \rangle := (u_h,v_h)_v$ for all
$u_h,v_h \in V_h$,
$\langle P_q \mathrm{p}, \mathrm{q} \rangle := (p_h,q_h)_q$ for all
$p_h,q_h \in Q_h$, and
$\langle B\mathrm{v}, \mathrm{q}\rangle := b_h(v_h, q_h)$ for all
$v_h \in V_h$ and $q_h \in Q_h$. Furthermore, we recall the following
result of \cite[Corollary 1]{baerland2020observation}. Let $P_{q,1}$
and $P_{q,0}$ be matrices defined for all $p_h,q_h \in Q_h$ by
$\langle P_{q,1}\mathrm{p},\mathrm{q}\rangle := (\nabla p_h, \nabla
q_h)$ and
$\langle P_{q,0} \mathrm{p}, \mathrm{q}\rangle := (p_h,q_h)$,
respectively. Then it can be shown that
$P_q^{-1} = \tau P_{q,1}^{-1} + \nu P_{q,0}^{-1}$.

\subsection{Preconditioning for stage-weighted Radau IIA RK
  time-stepping}

Let ${\bf u},{\bf v}$ be the expansion coefficients related to
$\boldsymbol{u}_h, \boldsymbol{v}_h \in \boldsymbol{V}_h$ and
${\bf p},{\bf q}$ be the expansion coefficients related to
$\boldsymbol{p}_h,\boldsymbol{q}_h \in \boldsymbol{Q}_h$. Let
$P_{\boldsymbol{v}}$ and $P_{\boldsymbol{q}}$ be matrices defined such
that
\begin{equation*}
  \langle P_{\boldsymbol{v}} {\bf u}, {\bf v} \rangle := (\boldsymbol{u}_h, \boldsymbol{v}_h)_{r,\boldsymbol{v}}
  \ \forall \boldsymbol{u}_h,\boldsymbol{v}_h \in \boldsymbol{V}_h,
  \qquad
  \langle P_{\boldsymbol{q}} {\bf p}, {\bf q} \rangle := (\boldsymbol{p}_h, \boldsymbol{q}_h)_{r,\boldsymbol{q}}
  \ \forall \boldsymbol{p}_h,\boldsymbol{q}_h \in \boldsymbol{Q}_h.
\end{equation*}
Furthermore, let $\boldsymbol{B}$ be the matrix such that
$\langle \boldsymbol{B}{\bf v}, {\bf q}\rangle :=
\boldsymbol{b}_r(\boldsymbol{v}_h, \boldsymbol{q}_h)$ for all
$\boldsymbol{v}_h \in \boldsymbol{V}_h$ and
$\boldsymbol{q} \in \boldsymbol{Q}_h$. By
\cref{thm:fovequiv,ss:radauiia}, the preconditioner defined by
\begin{equation}
  \label{eq:stokespreconradau}
  \mathcal{P}_U^{rk} :=
  \begin{bmatrix}
    \omega P_{\boldsymbol{v}} & \boldsymbol{B}^T
    \\
    0 & \rho P_{\boldsymbol{q}}
  \end{bmatrix},
\end{equation}
is FOV-equivalent with the sytem matrix associated with
\cref{eq:stokes-wradauiia}. We observe that $P_{\boldsymbol{v}}$,
$P_{\boldsymbol{q}}$, and $\boldsymbol{B}$ are block diagonal matrices
with one block per stage. Indeed, the $i^{th}$ block of these matrices
is given by $\hat{w}_iP_v$, $\hat{w}_iP_q$, and $\hat{w}_iB$,
respectively, for $i=1,\hdots,s$.

The inverse of $\mathcal{P}_U^{rk}$ is given by
\begin{equation*}
  (\mathcal{P}_U^{rk})^{-1} =
  \begin{bmatrix}
    \omega^{-1} P_{\boldsymbol{v}}^{-1} & -(\rho\omega)^{-1} P_{\boldsymbol{v}}^{-1} \boldsymbol{B}^T P_{\boldsymbol{q}}^{-1}
    \\
    0 & \rho^{-1} P_{\boldsymbol{q}}^{-1}
  \end{bmatrix}.
\end{equation*}
For the implementation of this preconditioner we then observe that
$P_{\boldsymbol{q}}^{-1}$ is a block diagonal matrix with one block
per stage and the $i^{th}$ block is given by
$\rho^{-1}\hat{w}_i^{-1}P_q^{-1} = \rho^{-1}\hat{w}_i^{-1}\tau
P_{q,1}^{-1} + \rho^{-1}\hat{w}_i^{-1}\nu P_{q,0}^{-1}$.

\begin{remark}
  Setting the weights $\hat{w}_i=1$ in \cref{eq:stokes-wradauiia} we
  obtain \cref{eq:stokes-radauiiaTH}. It can be shown that
  $\boldsymbol{\widetilde{a}}_r(\cdot, \cdot)$ is inf-sup stable but
  not necessarily coercive. In this case, one could replace
  $P_{\boldsymbol{v}}$ in \cref{eq:stokespreconradau} by the matrix
  associated with $\boldsymbol{\widetilde{a}}_r(\cdot, \cdot)$ and
  still have an FOV-equivalent preconditioner, see Remark
  \ref{rem:LogWathgeneral}. However, the matrix associated with
  $\boldsymbol{\widetilde{a}}_r(\cdot, \cdot)$ couples the stages. The
  proposed preconditioner \cref{eq:stokespreconradau} is block
  diagonal with one block per stage.
\end{remark}

\subsection{Preconditioning for weighted DG time-stepping}

Let ${\bf u},{\bf v}$ be the expansion coefficients related to
$\boldsymbol{u}_h, \boldsymbol{v}_h \in \mathcal{V}_h^n$ and
${\bf p},{\bf q}$ be the expansion coefficients related to
$\boldsymbol{p}_h, \boldsymbol{q}_h \in \mathcal{Q}_h^n$. Let
$P_{v,t}$ and $P_{q,t}$ be matrices defined such that
\begin{equation*}
  \langle P_{v,t} {\bf u}, {\bf v} \rangle := (u_h, v_h)_{v,t}
  \ \forall u_h,v_h \in \mathcal{V}_h^n,
  \qquad
  \langle P_{q,t} {\bf p}, {\bf q} \rangle := (p_h, q_h)_{q,t}
  \ \forall p_h,q_h \in \mathcal{Q}_h^n.
\end{equation*}
Furthermore, let $\boldsymbol{\hat{B}}$ be the matrix such that
$\langle \boldsymbol{\hat{B}}{\bf v}, {\bf q}\rangle :=
\boldsymbol{\hat{b}}(v_h,q_h)$ for all $v_h \in \mathcal{V}_h^n$ and
$q_h \in \mathcal{Q}_h$. By \cref{thm:fovequiv,sss:wspacetime}, the
preconditioner defined by
\begin{equation}
  \label{eq:ststokesprecon}
  \mathcal{P}_U^{st} :=
  \begin{bmatrix}
    \omega P_{v,t} & \boldsymbol{\hat{B}}^T
    \\
    0 & \rho P_{q,t}
  \end{bmatrix},
\end{equation}
is FOV-equivalent with the system matrix associated with
\cref{eq:discreteStokes}.

The inverse of $\mathcal{P}_U^{st}$ is given by
\begin{equation*}
  (\mathcal{P}_U^{st})^{-1} =
  \begin{bmatrix}
    \omega^{-1} P_{v,t}^{-1} & -(\rho\omega)^{-1} P_{v,t}^{-1} \boldsymbol{\hat{B}}^T P_{q,t}^{-1}
    \\
    0 & \rho^{-1} P_{q,t}^{-1}
  \end{bmatrix}.
\end{equation*}
For the implementation of (the action of) the inverse of
$\mathcal{P}_U^{st}$ we have the following useful characterization of
$P_{q,t}^{-1}$.

\begin{lemma}
  \label{lem:charac-st-preconditioner}
  Let $P_{qt,1}$ and $P_{qt,0}$ be matrices defined for all
  $p_h,q_h \in \mathcal{Q}_h$ by
  $\langle P_{qt,1}{\bf p}, {\bf q} \rangle := \int_{I_n} (\nabla p_h,
  \nabla q_h) \kappa_n \dif t$ and
  $\langle P_{qt,0}{\bf p}, {\bf q} \rangle := \int_{I_n} (p_h, q_h)
  \kappa_n \dif t$. Then
  \begin{equation*}
    P_{q,t}^{-1} = \tau P_{qt,1}^{-1} + \nu P_{qt,0}^{-1}.
  \end{equation*}
\end{lemma}
\begin{proof}
  We follow \cite[Corollary 1]{baerland2020observation}. For any
  $p_h\in \mathcal{Q}_h$, we have
  \begin{equation*}
    \langle p_h, p_h \rangle_{q,t}
    = \int_{I_n} \inf_{p_h = p_{h,1} + p_{h,0} }
    \del[1]{\tau^{-1}(\nabla p_{h,1}, \nabla p_{h,1}) 
      + \nu^{-1}(p_{h,0}, q_{h,0})} \kappa_n \dif{t}.
  \end{equation*}
  Applying \Cref{lem:inf-equality},
  \begin{equation*}
    \langle p_h, p_h \rangle_{q,t}
    = \inf_{ p_h = p_{h,1} + p_{h,0}}
    \del[3]{\tau^{-1}\int_{I_n} (\nabla p_{h,1}, \nabla p_{h,1}) \kappa_n \dif{t} 
      + \nu^{-1}\int_{I_n} (p_{h,0}, p_{h,0}) \kappa_n \dif{t}}.
  \end{equation*}
  Denoting the expansion coefficients related to
  $p_{h,1} \in \mathcal{Q}_h$ by ${\bf p}_1$ and the expansion
  coefficients related to $p_{h,0} \in \mathcal{Q}_h$ by ${\bf p}_0$,
  we write the above as
  \begin{equation*}
    \langle P_{q,t} {\bf p}, {\bf p}\rangle
    = \inf_{{\bf p} = {\bf p}_{1} + {\bf p}_{0}}\del[2]{
      \tau^{-1} \langle P_{qt,1} {\bf p}_{1}, {\bf p}_{1}\rangle
      + \nu^{-1} \langle P_{qt,0} {\bf p}_{0}, {\bf p}_{0}\rangle
    }.
  \end{equation*}
  The rest of the proof is analogous to \cite[Corollary
  1]{baerland2020observation}.
\end{proof}

\section{Numerical examples}
\label{s:numex}

For the numerical examples in this section we consider two- and
three-stage (stage-weighted) Radau IIA RK + Taylor--Hood
\cref{eq:stokes-wradauiia} and weighted space-time
\cref{eq:discreteStokes} discretizations of the time-dependent Stokes
problem. For the spatial discretization we consider $P^2$--$P^1$
Taylor--Hood element.

For the stage-weighted Radau IIA RK + Taylor--Hood discretization
\cref{eq:stokes-wradauiia} we require $C_r$ in
\cref{lem:coercivity-wradauiia} to be symmetric positive definite. We
choose $\widehat{W} = W$ with $W$ defined in \cref{lem:HHinvW}. It is
interesting to note that with this choice $C_r$ is symmetric positive
definite for $s=2,\hdots,7$. Indeed, let $\lambda_{\min}^{[s]}$ be the
smallest eigenvalue of $C_r$ of an $s$-stage Radau IIA RK method. Then
$\lambda^{[2]}_{\min} = 0.2603$, $\lambda^{[3]}_{\min} = 0.0842$,
$\lambda^{[4]}_{\min} = 0.0461$, $\lambda^{[5]}_{\min} = 0.0247$,
$\lambda^{[6]}_{\min} = 0.0121$, and $\lambda^{[7]}_{\min} = 0.0037$,
which are all positive. For $s=8,9,10$ we find
$\lambda^{[8]}_{\min} = -0.0023$, $\lambda^{[9]}_{\min} = -0.00768$,
and $\lambda^{[10]}_{\min} = -0.0105$ and so for $8 \le s \le 10$ one
would need to find a different matrix $\widehat{W}$. (The $s$-stage
Radau IIA RK tableaus for $s=2,\hdots,10$ were computed using
\cite{Kraus:2020:RungeKutta}.) For the 2-stage Radau IIA RK method we
therefore choose $\widehat{W} = \text{diag}(\hat{w}_1,\hat{w}_2)$ with
$\hat{w}_1 = (e-2)^{-1}e^{-1} \approx 0.51216588$ and
$\hat{w}_2 = (e-2)^{-1}(e+e^{-1}-3) \approx 0.11995468$, and for a
3-stage Radau IIA RK method we choose
$\widehat{W} = \text{diag}(\hat{w}_1,\hat{w}_2,\hat{w}_3)$ with
$\hat{w}_1=0.29272466$, $\hat{w}_2=0.29076228$, and
$\hat{w}_3=0.04863362$.

In what follows we will refer to \cref{eq:stokes-wradauiia}, in which
we use $s$-stages, as rk($s$) and we will refer to
\cref{eq:discreteStokes}, in which we use polynomials of degree $k_t$
in time, as dg($k_t$).
  
The simulations were performed using the finite element software
package \texttt{NGSolve} \cite{schoberl2014c++} together with
\texttt{NETGEN} \cite{schoberl1997netgen} to create the unstructured
meshes. For the dg($k_t$) implementation we use the NGSolve add-on
library \texttt{ngsxfem} \cite{ngsxfem}. We use the GMRES
implementation of NGSolve with a relative residual tolerance of
$10^{-8}$.

We consider two versions of the preconditioners $\mathcal{P}_U^{rk}$
(defined in \cref{eq:stokespreconradau}) and $\mathcal{P}_U^{st}$
(defined in \cref{eq:ststokesprecon}): (1) an `exact' version in which
a sparse direct solver is used to compute the inverses of each block;
and (2) an `inexact' version. The inexact version uses the Balancing
Domain Decomposition by Constraints (BDDC) method
\cite{dohrmann2003preconditioner} to approximate the action of the
inverse of the blocks associated with the `velocity'
degrees-of-freedom and to approximate the action of the inverse of the
blocks $P_{qt,1}$ (dg($k_t$)) and $P_{q,1}$ (rk($s$)). A direct solver
is used to compute the inverse of the mass matrices $P_{q,0}$ (rk($s$))
and $P_{qt,0}$ (dg($k_t$)) in two dimensions, while a Jacobi
preconditioner is used in three dimensions.

As test case we consider the backward facing step (BFS) problem
following the setup in \cite{fu2023uniform}. In two dimensions, we set
$\Omega := ([0.5,4] \times [0,0.5]) \cup ([0,4] \times [0.5,1])$ and
in three dimensions
$\Omega := (([0.5,4] \times [0,0.5]) \cup ([0,4] \times [0.5,1]))
\times [0,1]$. The boundary of the domains is divided into an inflow
boundary $\Gamma_{in}$ (at $x=0$), an outflow boundary $\Gamma_{out}$
(at $x=4$), and wall boundaries
$\Gamma_{wall} = \partial\Omega \backslash (\Gamma_{in} \cup
\Gamma_{out})$. On the inflow boundary, we impose
$u=[16(1-x_2)(x_2-0.5),0]^T$ in two dimensions and
$u=[64(1-x_2)(x_2-0.5)x_3(1-x_3), 0, 0]^T$ in three dimensions. On the
outflow boundary, we impose the Neumann boundary condition
$(-\nu \nabla u + pI)n = 0$, in which $n$ is the outward pointing
normal vector, and on the wall boundaries we impose $u=0$. Per test
case we perform only one time step. Since the outflow boundary is a
Neumann boundary condition we require to modify $(\cdot, \cdot)_q$ in
\cref{eq:innerprodvq} to
$(p_h, q_h)_q := \inf_{\substack{p_h = p_{h,1} + p_{h,2} \\ q_h =
    q_{h,1} + q_{h,2}}} \del[1]{\tau^{-1}(\nabla p_{h,1}, \nabla
  q_{h,1}) + \tau^{-1}h^{-1}(p_{h,1},q_{h,1})_{L^2(\Gamma_{out})} +
  \nu^{-1}(p_{h,2}, q_{h,2})}$.

\begin{figure}[tbp]
  \centering
  \includegraphics[width=1\linewidth]{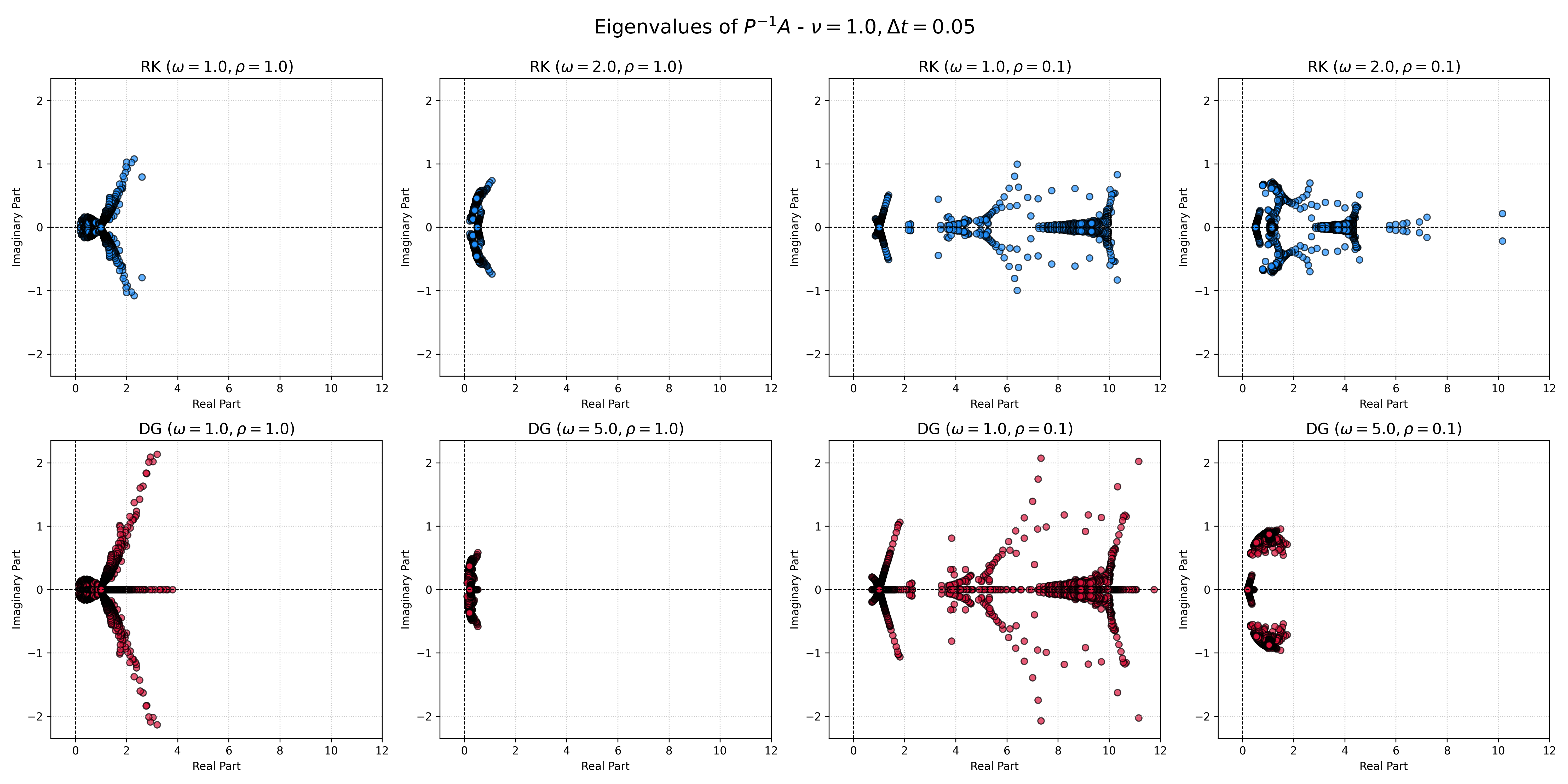}
  \caption{Two dimensional BFS problem. Distribution of the
    eigenvalues of $\mathcal{P}^{-1}\mathcal{A}$ for different
    $(\omega,\rho)$-pairs for rk($2$) and dg($2$). The eigenvalues are
    computed on a mesh consisting of 408 cells and with $\nu = 1$ and
    $\Delta t = 0.05$.}
  \label{fig:eigenvalueclustering} 
\end{figure}

We first study the performance of the `exact' form of the
preconditioners for different scaling parameters $\rho$ and
$\omega$. In \cref{tab:scaling-omega-rho-dg-rk} we list for the two
dimensional problem the number of GMRES iterations to reach the
desired tolerance for one time step. We do this for rk($2$), rk($3$),
dg($2$), and dg($3$). We fix the number of cells in the mesh to 2166
and consider $\rho \in \cbr[0]{0.01, 0.05, 0.1, 0.5, 1.0}$,
$\omega \in \cbr[0]{1, 2, 5, 10, 15}$,
$\Delta t \in \cbr[0]{10^{-1}, 10^{-2}, 10^{-3}}$, and
$\nu \in \cbr[0]{1, 10^{-3}}$. For the majority of cases we observe
that the lowest iteration count is obtained with
$\omega \in \cbr[0]{2,5}$ and $\rho \in \cbr[0]{0.1, 0.5}$. From
\cref{thm:fovequiv} we observe that large $\omega$ has the effect of
reducing the difference between $c_{\star}$ and $\bar{c}_{\star}$,
which implies a better clustering of the eigenvalues of
$\mathcal{P}_U^{-1}\mathcal{A}$. This is also demonstrated for the two
dimensional BFS problem in \cref{fig:eigenvalueclustering} for rk($2$)
and dg($2$). However, increasing $\omega$ also has the effect of
shifting the eigenvalues closer to the $y$-axis. The eigenvalues are
shifted away from the $y$-axis by decreasing $\rho$.

\begin{table}[tbp]
  \centering
  \resizebox{\textwidth}{!}{%
    \begin{tabular}{c|c|c|ccccc|ccccc|ccccc|ccccc|}
      \hline
      & & & \multicolumn{5}{c|}{rk(2)} & \multicolumn{5}{c|}{rk(3)} & \multicolumn{5}{c|}{dg(2)} & \multicolumn{5}{c|}{dg(3)} \\
      \cline{4-23}
      & & & \multicolumn{20}{c|}{$\omega$} \\
      \cline{4-23}
      $\nu$ & $\Delta t$ & $\rho$ & $1$ & $2$ & $5$ & $10$ & $15$ & $1$ & $2$ & $5$ & $10$ & $15$ & $1$ & $2$ & $5$ & $10$ & $15$ & $1$ & $2$ & $5$ & $10$ & $15$ \\
      \hline
      $1$ & $10^{-1}$ & $0.01$ &   54 & 54 & 58 & 60 & 62 & 65 & 67 & 75 & 77 & 78 & 68 & 68 & 76 & 80 & 81 & 81 & 79 & 87 & 92 & 96  \\
      & & $0.05$ & 41 & 47 & 52 & 59 & 62 & 50 & 50 & 61 & 69 & 76 & 52 & 51 & 61 & 70 & 78 & 62 & 56 & 68 & 78 & 87 \\
      & & $0.1$ & {\bf 38} & 42 & 51 & 60 & 64 & 47 & {\bf 46} & 56 & 69 & 76 & 48 & {\bf 47} & 57 & 69 & 77 & 60 & {\bf 52} & 61 & 78 & 87 \\
      & & $0.5$ &  39 & 39 & 52 & 59 & 63 & 58 & 51 & 57 & 67 & 73 & 59 & 52 & 58 & 68 & 74 & 86 & 73 & 68 & 76 & 84 \\
      & & $1.0$ &  46 & 45 & 52 & 58 & 61 & 76 & 65 & 61 & 66 & 69 & 77 & 67 & 62 & 68 & 71 & 113 & 96 & 82 & 78 & 82 \\
      \cline{2-23}
      & $10^{-2}$ & $0.01$ &  65 & 60 & 63 & 68 & 68 & 83 & 77 & 77 & 81 & 82 & 84 & 77 & 77 & 80 & 82 & 110 & 95 & 91 & 96 & 100 \\
      & & $0.05$ &  48 & 46 & 51 & 56 & 63 & 59 & 50 & 58 & 68 & 75 &  59 & 50 & 59 & 68 & 75 & 76 & 61 & 64 & 74 & 85 \\
      & & $0.1$ &  42 & {\bf 41} & 49 & 59 & 64 & 54 & {\bf 44} & 52 & 65 & 75 & 54 & {\bf 45} & 52 & 65 & 75 & 75 & {\bf 58} & 61 & 74 & 84 \\
      & & $0.5$ & 45 & 40 & 50 & 60 & 63 & 79 & 61 & 59 & 65 & 71 & 80 & 64 & 61 & 67 & 72 & 130 & 98 & 81 & 81 & 83 \\
      & & $1.0$ & 57 & 47 & 51 & 60 & 65 & 109 & 85 & 69 & 69 & 71 & 109 & 87 & 72 & 71 & 74 & 183 & 140 & 105 & 93 & 91 \\
      \cline{2-23}
      & $10^{-3}$ & $0.01$ &  54 & 49 & 48 & 51 & 50 & 79 & 68 & 65 & 67 & 68 & 80 & 70 & 66 & 67 & 68 & 107 & 88 & 80 & 82 & 83 \\
      & & $0.05$ &  42 & 37 & 43 & 46 & 50 & 59 & 48 & 50 & 56 & 62 &  59 & 48 & 51 & 57 & 63 & 80 & 61 & 58 & 64 & 73 \\
      & & $0.1$ &   37 & {\bf 36} & 39 & 47 & 53 & 56 & {\bf 45} & 47 & 57 & 63 & 56 & {\bf 47} & {\bf 47} & 57 & 65 & 79 & 61 & {\bf 56} & 65 & 73 \\
      & & $0.5$ &  43 & {\bf 36} & 42 & 48 & 53 & 82 & 60 & 54 & 58 & 63 & 82 & 61 & 57 & 60 & 67 & 140 & 101 & 78 & 73 & 75 \\
      & & $1.0$ & 54 & 44 & 45 & 50 & 53 & 115 & 85 & 66 & 63 & 65 & 113 & 86 & 67 & 66 & 66 & 204 & 149 & 104 & 89 & 85 \\
      \hline
      $10^{-3}$ & $10^{-1}$ & $0.01$ &  93 & 92 & 95 & 97 & 100 & 138 & 126 & 124 & 131 & 136 & 136 & 123 & 120 & 128 & 133 & 183 & 159 & 152 & 158 & 165 \\
      & & $0.05$ &  74 & 70 & 79 & 85 & 92 & 97 & 84 & 86 & 98 & 109 & 94 & 80 & 84 & 95 & 105 & 124 & 100 & 97 & 110 & 120 \\
      & & $0.1$ &  63 & 60 & 70 & 85 & 98 & 84 & {\bf 68} & 73 & 89 & 107 & 84 & 66 & 71 & 86 & 103 & 116 & 88 & {\bf 82} & 96 & 114 \\
      & & $0.5$ &  63 & {\bf 52} & 64 & 83 & 96 & 106 & 76 & 65 & 78 & 91 & 105 & 75 & {\bf 63} & 76 & 89 & 172 & 122 & 85 & 87 & 97 \\
      & & $1.0$ &  71 & 58 & 67 & 81 & 90 & 139 & 101 & 76 & 78 & 84 & 138 & 100 & 75 & 77 & 84 & 243 & 173 & 115 & 98 & 98 \\
      \cline{2-23}
      & $10^{-2}$ & $0.01$ &  95 & 85 & 86 & 88 & 89 & 130 & 109 & 100 & 105 & 109 & 126 & 107 & 96 & 101 & 105 & 170 & 142 & 115 & 118 & 122 \\
      & & $0.05$ &   79 & 72 & 74 & 79 & 81 & 100 & 81 & 73 & 80 & 86 & 102 & 79 & 69 & 77 & 82 & 136 & 103 & 77 & 82 & 90 \\
      & & $0.1$ & 74 & {\bf 62} & 68 & 75 & 82 & 98 & 75 & {\bf 61} & 70 & 78 & 98 & 75 & {\bf 58} & 67 & 75 & 139 & 106 & 75 & {\bf 71} & 77 \\
      & & $0.5$ & 78 & 64 & {\bf 62} & 78 & 88 & 127 & 95 & 72 & 73 & 77 & 125 & 94 & 69 & 71 & 75 & 204 & 151 & 101 & 89 & 90 \\
      & & $1.0$ &  86 & 69 & 71 & 81 & 87 & 160 & 118 & 87 & 84 & 84 & 156 & 116 & 84 & 82 & 83 & 276 & 203 & 133 & 112 & 108 \\
      \cline{2-23}
      & $10^{-3}$ & $0.01$ & 75 & 67 & 59 & 63 & 62 & 106 & 86 & 75 & 77 & 81 & 110 & 86 & 72 & 76 & 80 & 133 & 114 & 90 & 87 & 91 \\
      & & $0.05$ & 65 & 53 & 54 & 57 & 59 & 89 & 72 & 59 & 63 & 68 & 91 & 72 & 58 & 61 & 67 & 123 & 97 & 71 & 67 & 71 \\
      & & $0.1$ &  62 & {\bf 50} & {\bf 50} & 56 & 57 & 93 & 73 & 59 & {\bf 57} & 62 & 95 & 73 & 60 & {\bf 57} & 61 & 135 & 106 & 76 & 71 & {\bf 68} \\
      & & $0.5$ & 67 & 54 & {\bf 50} & 55 & 58 & 125 & 97 & 72 & 70 & 70 & 124 & 97 & 72 & 70 & 71 & 200 & 156 & 110 & 93 & 92 \\
      & & $1.0$ & 75 & 59 & 57 & 59 & 62 & 157 & 122 & 90 & 84 & 83 & 154 & 120 & 88 & 83 & 83 & 261 & 203 & 142 & 119 & 113 \\
      \hline
    \end{tabular} 
  }
  \caption{Two dimensional BFS problem on a mesh with 2166
    cells. Number of GMRES iteration needed to reach the desired
    tolerance when using the `exact' form of the preconditioner with
    varying scaling parameters $\rho$ and $\omega$.}
  \label{tab:scaling-omega-rho-dg-rk}    
\end{table}

We now consider $h$-robustness of the preconditioners for rk($2$),
rk($3$), dg($2$), and dg($3$). We fix $\nu = 1$ and
$\Delta t = 10^{-3}$. In \cref{tab:h-robustness} we list for the two
and three dimensional BFS problems the number of GMRES iterations
needed to reach the desired tolerance. We observe that, when using the
exact versions of the preconditioners, the iteration count decreases
at the number of mesh cells increases. When using the inexact versions
of the preconditioners, the variation in number of iterations is very
small as the number of cells in the mesh changes.

\begin{table}[tpb]
  \centering
  \resizebox{\textwidth}{!}{%
    \begin{tabular}{c | c c c c c }
      \hline
      \multicolumn{6}{c}{2D}\\
      \hline
      $ \text{scheme} \setminus \text{cells}$ & $102$ & $408$ & $1632$ & $6528$ & $26112$  \\
      \hline
      rk(2) ($\omega = 2$, $\rho=0.1$) & 35 (36) & 37 (37) & 38 (42) & 34 (46) & 31 (49) \\
      rk(3) ($\omega = 2$, $\rho=0.1$) & 56 (57) & 49 (52) & 48 (54) & 43 (57) & 40 (62) \\
      dg(2) ($\omega = 5$, $\rho=0.1$) & 46 (45) & 50 (52) & 52 (58) & 47 (62) & 42 (61) \\
      dg(3) ($\omega = 2$, $\rho=0.1$) & 60 (61) & 60 (61) & 61 (67) &  54 (66) & 49 (68) \\ 
      \hline
      \multicolumn{6}{c}{3D}\\
      \hline
      $ \text{scheme} \setminus \text{cells}$ & $37$ & $296$ & $2368$ & $18944$ & $151552$  \\
      \hline
      rk(2) ($\omega = 2$, $\rho=0.1$)   & 140 (140) & 60 (62) & 60 (65) & 58 (68) & 51 (63) \\
      rk(3) ($\omega = 2$, $\rho=0.1$)   & 208 (211) & 104 (112) & 93 (104) & 81 (96) & 62 (76) \\
      dg(2) ($\omega = 5$, $\rho=0.05$)  & 206 (209) & 71 (72) & 82 (88) & 89 (96) & 77 (100)\\
      dg(3) ($\omega = 10$, $\rho=0.05$) & 277 (280) & 89 (95) & 105 (111) & 110 (125) & 97 (110) \\ 
      \hline
    \end{tabular}}
  \caption{The two and three dimensional BFS problems. Number of GMRES
    iterations needed to reach the desired tolerance when using the
    exact and inexact (in parentheses) versions of the
    preconditioners. We set $\Delta t = 10^{-3}$, $\nu = 1$, and vary
    the number of cells in the mesh.}
  \label{tab:h-robustness}
\end{table}

Finally, we demonstrate robustness of the preconditioners with respect
to the time step $\Delta t$ and viscosity $\nu$ for rk($2$), rk($3$),
dg($2$), and dg($3$) in \cref{tab:BFS_dt_vs_nu}. In two dimensions we
consider a mesh consisting of 8724 cells and in three dimensions a
mesh consisting of 4693 cells. We observe some dependence of the
solver on $\Delta t$ but given the variation in $\Delta t$ (from
$\Delta t = 10^{-4}$ to $\Delta t = 1$) this variation is very
small. For a given $\Delta t$ there is also some variation in
iteration count as viscosity decreases from $\nu=1$ to $\nu=10^{-4}$
but again, this variation is small given the large variation in
$\nu$. For viscosities between $\nu = 10^{-8}$ and $\nu = 10^{-4}$
there is no change in iteration count.

\begin{table}[tbp]
  \centering
  \resizebox{\textwidth}{!}{%
    \begin{tabular}{l | ccccc | ccccc}
      \hline
      & \multicolumn{5}{c|}{rk($2$), $\omega= 2$, $\rho= 0.1$} & \multicolumn{5}{c}{rk($3$), $\omega= 2$, $\rho= 0.1$} \\
      \cline{2-6} \cline{7-11}
      $\nu \setminus \Delta t $ & $10^{0}$ & $10^{-1}$ & $10^{-2}$ & $10^{-3}$ & $10^{-4}$ & $10^{0}$ & $10^{-1}$ & $10^{-2}$ & $10^{-3}$ & $10^{-4}$ \\
      \hline
      \multicolumn{11}{c}{2D} \\
      \hline
      $1$       & 49 (52) & 40 (49) & 37 (42) & 36 (36) & 38 (39) & 52 (59) & 46 (53) & 44 (52) & 45 (46) & 49 (50) \\
      $10^{-2}$ & 43 (48) & 41 (43) & 45 (45) & 47 (47) & 46 (52) & 52 (60) & 54 (55) & 58 (59) & 68 (70) & 71 (78) \\
      $10^{-4}$ & 58 (59) & 59 (60) & 54 (60) & 50 (56) & 48 (53) &  77 (78) & 86 (89) & 80 (91) & 73 (85) & 67 (81) \\
      $10^{-6}$ & 68 (76) & 62 (71) & 54 (62) & 50 (57) & 48 (53) &  100 (115) & 89 (107) & 78 (94) & 71 (85) & 67 (81) \\
      $10^{-8}$ & 67 (78) & 60 (71) & 54 (63) & 50 (57) & 48 (53) & 94 (119) & 88 (107) & 76 (94) & 71 (85) & 67 (81) \\
      \hline
      \hline
      & \multicolumn{5}{c|}{rk($2$), $\omega= 2$, $\rho= 0.1$} & \multicolumn{5}{c}{rk($3$), $\omega= 2$, $\rho= 0.1$} \\
      \hline
      $\nu \setminus \Delta t $ & $10^{0}$ & $10^{-1}$ & $10^{-2}$ & $10^{-3}$ & $10^{-4}$ & $10^{0}$ & $10^{-1}$ & $10^{-2}$ & $10^{-3}$ & $10^{-4}$ \\
      \hline
      \multicolumn{11}{c}{3D} \\
      \hline
      $1$       & 60 (77) & 53 (68) & 47 (57) & 48 (52) & 52 (53) & 71 (86) & 63 (79) & 68 (83) & 100 (109) & 146 (149) \\
      $10^{-2}$ & 54 (62) & 51 (56) & 54 (55) & 55 (55) & 55 (55) & 74 (91) & 106 (116) & 151 (153) & 156 (157) & 148 (148) \\
      $10^{-4}$ & 67 (68) & 62 (62) & 56 (56) & 55 (55) & 55 (55) & 166 (169) & 167 (168) & 160 (160) & 151 (151) & 147 (147) \\
      $10^{-6}$ & 65 (65) & 62 (62) & 56 (56) & 55 (55) & 55 (55) & 168 (168) & 165 (165) & 159 (159) & 149 (149) & 147 (147) \\
      $10^{-8}$ & 65 (65) & 62 (62) & 56 (56) & 55 (55) & 55 (55) & 166 (166) & 165 (165) & 158 (158) & 149 (149) & 147 (147) \\
      \hline
      \hline
      & \multicolumn{5}{c|}{dg($2$), $\omega = 5$, $\rho = 0.1$ } & \multicolumn{5}{c}{dg($3$), $\omega = 5$, $\rho = 0.1$} \\
      \cline{2-6} \cline{7-11}
      $\nu \setminus \Delta t $ & $10^{0}$ & $10^{-1}$ & $10^{-2}$ & $10^{-3}$ & $10^{-4}$ & $10^{0}$ & $10^{-1}$ & $10^{-2}$ & $10^{-3}$ & $10^{-4}$ \\
      \hline
      \multicolumn{11}{c}{2D} \\
      \hline
      $1$       & 62 (62) & 57 (57) & 52 (52) & 47 (47) & 51 (51) &  66 (66) & 61 (61) & 61 (61) & 56 (56) & 59 (59)  \\
      $10^{-2}$ & 63 (63) & 62 (62) & 62 (62) & 55 (55) & 59 (59) & 73 (73) & 71 (71) & 72 (72) & 71 (71) & 75 (75) \\
      $10^{-4}$ & 81 (81) & 67 (67) & 63 (63) & 61 (61) & 60 (60) &93 (93) & 85 (85) & 82 (82) & 76 (76) & 76 (76) \\
      $10^{-6}$ & 77 (77) & 70 (70) & 64 (64) & 61 (61) & 60 (60) & 100 (100) & 93 (93) & 82 (82) & 76 (76) & 76 (76) \\
      $10^{-8}$ & 77 (77) & 70 (70) & 64 (64) & 61 (61) & 60 (60) & 100 (100) & 92 (92) & 82 (82) & 76 (76) & 76 (76) \\
      \hline
      \hline
      & \multicolumn{5}{c|}{dg($2$), $\omega = 5$, $\rho = 0.05$} & \multicolumn{5}{c}{dg($3$), $\omega = 10$, $\rho = 0.05$} \\
      \hline
      $\nu \setminus \Delta t $ & $10^{0}$ & $10^{-1}$ & $10^{-2}$ & $10^{-3}$ & $10^{-4}$ & $10^{0}$ & $10^{-1}$ & $10^{-2}$ & $10^{-3}$ & $10^{-4}$ \\
      \hline
      \multicolumn{11}{c}{3D} \\
      \hline
      $1$       & 97 (139) & 90 (130) & 87 (107) & 87 (86) & 102 (103) & 115 (160) & 114 (158) & 112 (130) & 107 (112) & 122 (124)  \\
      $10^{-2}$ &  104 (124) & 98 (99) & 102 (104) & 108 (108) & 109 (109) & 139 (161) & 125 (130) & 122 (124) & 128 (128) & 130 (130) \\
      $10^{-4}$ & 115 (118) & 114 (114) & 109 (109) & 109 (109) & 109 (109) & 137 (139) & 134 (134) & 130 (130) & 130 (130) & 130 (130) \\
      $10^{-6}$ &  123 (123) & 114 (114) & 109 (109) & 109 (109) & 109 (109) & 145 (145) & 135 (135) & 130 (130) & 130 (130) & 130 (130) \\
      $10^{-8}$ &  123 (123) & 114 (114) & 109 (109) & 109 (109) & 109 (109) & 145 (145) & 135 (135) & 130 (130) & 130 (130) & 130 (130) \\
      \hline      
    \end{tabular}}
  \caption{The two and three dimensional BFS problems. Number of GMRES
    iterations needed to reach the desired tolerance when using the
    exact and inexact (in parentheses) versions of the
    preconditioners. This table verifies robustness of the
    preconditioners with respect to the time step and the
    viscosity. The spatial mesh is fixed with 8724 cells for the two
    dimensional runs and 4693 for the three dimensional runs.}
  \label{tab:BFS_dt_vs_nu}    
\end{table}

\section*{Acknowledgments}

SR acknowledges support from the Natural Sciences and Engineering
Research Council of Canada through the Discovery Grant program
(RGPIN-2023-03237).

\section*{Declarations}

ChatGPT 6.1 Sol was used to obtain a first draft of the mathematical
results in \cref{sss:wspacetimeRK} which were then checked and revised
by the authors. Microsoft Copilot Chat (Basic) was used in the
remainder of the manuscript to confirm and simplify some mathematical
results originally derived without AI assistance, and for some
mathematical exploration which was further worked out by the
authors. Gemini 1.5 Pro was used to improve codes used in this
paper. The authors assume responsibility for all content.

\bibliographystyle{plain}
\bibliography{references-arxiv}
\appendix
\section{Proof of \cref{thm:fovequiv}}
\label{ap:proofFOVequiv}

Before proving \cref{thm:fovequiv} we require the following results.

\begin{lemma}
  \label{lem:boundAinv}
  The $\alpha$ and $c_a$ be the constants from
  \cref{eq:aconditions}. Then,
  \begin{equation*}
    \alpha c_a^{-2} \langle P_v^{-1} {\bf u}, {\bf u} \rangle
    \le \langle A^{-1} {\bf u}, {\bf u} \rangle
    \le c_a \alpha^{-2} \langle P_v^{-1} {\bf u}, {\bf u} \rangle
    \qquad \forall {\bf u} \in \mathbb{R}^{n_v}.
  \end{equation*}
\end{lemma}
\begin{proof}
  Let ${\bf v} \in \mathbb{R}^{n_v}$ and ${\bf u} = A{\bf v}$. By
  \cref{eq:abs-coercivity} and Cauchy--Schwarz,
  \begin{equation}
    \label{eq:alphaPvvv}
    \alpha \langle P_v{\bf v}, {\bf v} \rangle
    \le \langle A{\bf v}, {\bf v} \rangle
    = \langle {\bf u}, {\bf v} \rangle
    \le \langle P_v^{-1}{\bf u}, {\bf u} \rangle^{1/2} \langle P_v{\bf v}, {\bf v} \rangle^{1/2}.
  \end{equation}
  Combining \cref{eq:a-boundedness} with \cref{eq:alphaPvvv},
  \begin{equation*}
    \langle A^{-1}{\bf u}, {\bf u} \rangle
    = \langle A {\bf v}, {\bf v} \rangle
    \le c_a \langle P_v{\bf v}, {\bf v} \rangle
    \le c_a \alpha^{-2} \langle P_v^{-1}{\bf u}, {\bf u} \rangle,
  \end{equation*}
  proving the upper bound. For the lower bound, using the definition
  of the dual norm and \cref{eq:a-boundedness},
  \begin{equation*}
    \norm[0]{{\bf u}}_{P_v^{-1}}
    = \sup_{\substack{{\bf w} \in \mathbb{R}^{n_v} \\ {\bf w} \ne {\bf 0}}} \frac{\langle {\bf u}, {\bf w} \rangle}{\norm[0]{{\bf w}}_{P_v}}
    = \sup_{\substack{{\bf w} \in \mathbb{R}^{n_v} \\ {\bf w} \ne {\bf 0}}} \frac{\langle A{\bf v}, {\bf w} \rangle}{\norm[0]{{\bf w}}_{P_v}}
    \le \sup_{\substack{{\bf w} \in \mathbb{R}^{n_v} \\ {\bf w} \ne {\bf 0}}} \frac{c_a \norm[0]{{\bf v}}_{P_v} \norm[0]{{\bf w}}_{P_v} }{\norm[0]{{\bf w}}_{P_v}},
  \end{equation*}
  so that
  $\langle P_v^{-1} {\bf u}, {\bf u} \rangle \le c_a^2 \langle P_v
  {\bf v}, {\bf v} \rangle$. Then, again using
  \cref{eq:abs-coercivity},
  \begin{equation*}
    \alpha c_a^{-2} \langle P_v^{-1} {\bf u}, {\bf u} \rangle
    \le \alpha \langle P_v{\bf v}, {\bf v} \rangle
    \le \langle A{\bf v}, {\bf v} \rangle
    = \langle A^{-1}{\bf w}, {\bf w} \rangle,
  \end{equation*}
  proving the lower bound.
\end{proof}

\begin{lemma}
  \label{lem:boundS}
  Let $S = BA^{-1}B^T$ and let $c_a,c_b,\alpha,\beta$ be the constants
  from \cref{eq:bconditions,eq:aconditions}. Then,
  \begin{equation*}
    \alpha c_a^{-2} \beta^2 \langle P_q {\bf q}, {\bf q} \rangle
    \le \langle S {\bf q}, {\bf q} \rangle
    \le c_a \alpha^{-2} c_b^2 \langle P_q{\bf q}, {\bf q} \rangle
    \qquad \forall {\bf q} \in \mathbb{R}^{n_q}.
  \end{equation*}
\end{lemma}
\begin{proof}
  Using \cref{lem:boundAinv},
  \begin{equation}
    \label{eq:Sqqbound}
    \langle S{\bf q}, {\bf q} \rangle
    = \langle BA^{-1}B^T{\bf q}, {\bf q} \rangle
    = \langle A^{-1}B^T{\bf q}, B^T{\bf q} \rangle
    \ge \alpha c_a^{-2} \langle P_v^{-1} B^T{\bf q}, B^T{\bf q} \rangle.
  \end{equation}
  Using \cref{eq:abs-infsup} and the definition of the dual norm, we
  note that
  $\norm[0]{B^T{\bf q}}_{P_v^{-1}} \ge \beta \norm[0]{{\bf
      q}}_{P_q}$. Combining with \cref{eq:Sqqbound},
  \begin{equation*}
    \langle S{\bf q}, {\bf q} \rangle
    \ge \alpha c_a^{-2} \beta^2 \langle P_q {\bf q}, {\bf q} \rangle,
  \end{equation*}
  proving the lower bound. For the upper bound, using again
  \cref{lem:boundAinv},
  \begin{equation}
    \label{eq:Sqqboundup}
    \langle S{\bf q}, {\bf q} \rangle
    = \langle BA^{-1}B^T{\bf q}, {\bf q} \rangle
    = \langle A^{-1}B^T{\bf q}, B^T{\bf q} \rangle
    \le c_a \alpha^{-2} \langle P_v^{-1} B^T{\bf q}, B^T{\bf q} \rangle.
  \end{equation}
  Using \cref{eq:b-boundedness} and the definition of the dual norm,
  $\norm[0]{B^T{\bf q}}_{P_v^{-1}} \le c_b \norm[0]{{\bf
      q}}_{P_q}$. Combining with \cref{eq:Sqqboundup},
  \begin{equation*}
    \langle S{\bf q}, {\bf q} \rangle
    \le c_a \alpha^{-2} c_b^2 \langle P_q{\bf q}, {\bf q} \rangle,
  \end{equation*}
  proving the upper bound.
\end{proof}

We now prove the ($\mathcal{P}^{-1}$-)FOV-equivalence of $\mathcal{A}$
and $\mathcal{P}_U$.

\begin{proof}[Proof of \cref{thm:fovequiv}]
  {\bf Step 1.} We start with some useful inequalities. Let
  $A_s = \tfrac{1}{2}(A + A^T)$ be the symmetric part of $A$ and note
  that
  $\langle A {\bf v}, {\bf v} \rangle = \langle A_s{\bf v}, {\bf v}
  \rangle$.  Let ${\bf v} = A_s^{-1}{\bf w}$ for any
  ${\bf v} \in \mathbb{R}^{n_v}$ and note that
  $\langle A{\bf v}, {\bf v} \rangle = \langle A_s{\bf v}, {\bf v}
  \rangle = \langle A_s^{-1}{\bf w}, {\bf w}\rangle$. Then, using
  Cauchy--Schwarz and \cref{eq:abs-coercivity},
  $ \langle A_s^{-1}{\bf w}, {\bf w} \rangle = \langle {\bf w}, {\bf
    v} \rangle \le \norm[0]{{\bf w}}_{P_v^{-1}} \norm[0]{{\bf
      v}}_{P_v} \le \norm[0]{{\bf w}}_{P_v^{-1}} \alpha^{-1/2}\langle
  A_s{\bf v}, {\bf v} \rangle^{1/2}$,
  so that
  \begin{equation}
    \label{eq:Asmin1wwbound}
    \langle A_s^{-1}{\bf w}, {\bf w} \rangle \le \alpha^{-1} \langle
    P_v^{-1}{\bf w}, {\bf w} \rangle. 
  \end{equation}
  By definition of the dual norm and \cref{eq:a-boundedness},
  \begin{equation}
    \label{eq:Avdualbound}
    \norm[0]{A{\bf w}}_{P_v^{-1}}
    = \sup_{\substack{{\bf v} \in \mathbb{R}^{n_v} \\ {\bf v} \ne {\bf 0}}}
    \frac{\langle A{\bf w}, {\bf v} \rangle}{\norm[0]{{\bf v}}_{P_v}}
    \le \sup_{\substack{{\bf v} \in \mathbb{R}^{n_v} \\ {\bf v} \ne {\bf 0}}}
    \frac{c_a \norm[0]{{\bf w}}_{P_v} \norm[0]{{\bf v}}_{P_v}}{\norm[0]{{\bf v}}_{P_v}}
    = c_a \norm[0]{{\bf w}}_{P_v}.
  \end{equation}
  Likewise, by definition of the dual norm and
  \cref{eq:b-boundedness},
  \begin{equation}
    \label{eq:Bvdualbound}
    \norm[0]{B{\bf v}}_{P_q^{-1}}
    = \sup_{\substack{{\bf q} \in \mathbb{R}^{n_q} \\ {\bf q} \ne {\bf 0}}}
    \frac{\langle B{\bf v}, {\bf q} \rangle}{\norm[0]{{\bf q}}_{P_q}}
    \le \sup_{\substack{{\bf q} \in \mathbb{R}^{n_q} \\ {\bf q} \ne {\bf 0}}}
    \frac{c_b \norm[0]{{\bf v}}_{P_v} \norm[0]{{\bf q}}_{P_q}}{\norm[0]{{\bf q}}_{P_q}}
    = c_b \norm[0]{{\bf v}}_{P_v}.
  \end{equation}
  {\bf Step 2.} We determine $c_{\star}$. Using
  \cref{eq:abs-coercivity,lem:boundAinv,lem:boundS}, and taking the
  infimum over non-zero vectors,
  \begin{equation}
    \label{eq:infcstar}
    \begin{split}
      &\inf_{{\bf x} \in \mathbb{R}^n} \frac{\langle AP_U^{-1}{\bf x}, {\bf x} \rangle_{\mathcal{P}^{-1}}}{\langle {\bf x}, {\bf x} \rangle_{\mathcal{P}^{-1}}}
      \\
      =&\inf_{({\bf v},{\bf q}) \in \mathbb{R}^{n_v} \times \mathbb{R}^{n_q}}
         \frac{\omega^{-2}\langle P_v^{-1}AP_v^{-1}{\bf v}, {\bf v} \rangle
         - \rho^{-1}\omega^{-2}\langle P_v^{-1}AP_v^{-1}B^TP_q^{-1}{\bf q}, {\bf v} \rangle
         + \rho^{-2}\omega^{-1}\langle P_q^{-1}BP_v^{-1}B^TP_q^{-1}{\bf q}, {\bf q} \rangle
         }{\omega^{-1} \langle P_v^{-1}{\bf v}, {\bf v} \rangle + \rho^{-1}\langle P_q^{-1} {\bf q}, {\bf q} \rangle }
      \\
      \ge&\inf_{({\bf v},{\bf q}) \in \mathbb{R}^{n_v} \times \mathbb{R}^{n_q}}
           \frac{\omega^{-2}\langle P_v^{-1}A_sP_v^{-1}{\bf v}, {\bf v} \rangle
           - \rho^{-1}\omega^{-2}\langle P_v^{-1}AP_v^{-1}B^TP_q^{-1}{\bf q}, {\bf v} \rangle
           + \rho^{-2}\omega^{-1}\langle P_q^{-1}BP_v^{-1}B^TP_q^{-1}{\bf q}, {\bf q} \rangle
           }{\alpha^{-1} \omega^{-1} \langle P_v^{-1}A_sP_v^{-1}{\bf v}, {\bf v} \rangle + \rho^{-1}\alpha^{-3}c_a^3 \beta^{-2}\langle P_q^{-1} B P_v^{-1}B^TP_q^{-1} {\bf q}, {\bf q} \rangle }
      \\
      \ge & \frac{1}{\tilde{c}}\inf_{({\bf v},{\bf q}) \in \mathbb{R}^{n_v} \times \mathbb{R}^{n_q}}
            \frac{\rho^2 \langle P_v^{-1}A_sP_v^{-1}{\bf v}, {\bf v} \rangle
            - \rho \langle P_v^{-1}AP_v^{-1}B^TP_q^{-1}{\bf q}, {\bf v} \rangle
            + \rho \omega\langle P_q^{-1}BP_v^{-1}B^TP_q^{-1}{\bf q}, {\bf q} \rangle
            }{\rho^2\langle P_v^{-1}A_sP_v^{-1}{\bf v}, {\bf v} \rangle + \rho \omega\langle P_q^{-1} B P_v^{-1}B^TP_q^{-1} {\bf q}, {\bf q} \rangle },        
    \end{split}
  \end{equation}
  where
  $\tilde{c} = \omega \alpha^{-1} \max\{1, \alpha^{-2}c_a^3\beta^{-2}
  \}$. To find a lower bound for the right hand side of
  \cref{eq:infcstar}, we follow the proof of \cite[Theorem
  4.3]{he2025preconditioner}, i.e., we find lower bounds for the
  generalized eigenvalues of
  \begin{equation*}
    \begin{bmatrix}
      \rho^2 P_v^{-1}A_sP_v^{-1}
      & -\tfrac{1}{2} \rho P_v^{-1}AP_v^{-1}B^TP_q^{-1}
      \\
      -\tfrac{1}{2} \rho P_q^{-1}B P_v^{-1}A^TP_v^{-1}
      &
        \rho \omega P_q^{-1}BP_v^{-1}B^TP_q^{-1}
    \end{bmatrix}
    \begin{bmatrix}
      {\bf v} \\ {\bf q}
    \end{bmatrix}
    =
    \lambda
    \begin{bmatrix}
      \rho^2 P_v^{-1}A_sP_v^{-1} & 0
      \\
      0 & \rho \omega P_q^{-1}BP_v^{-1}B^TP_q^{-1}
    \end{bmatrix}
    \begin{bmatrix}
      {\bf v} \\ {\bf q}
    \end{bmatrix}.
  \end{equation*}
  After algebraic manupulation, we find that $\lambda$ satisfies
  \begin{equation}
    \label{eq:lambdamin1}
    (\lambda - 1)^2 = \frac{\langle P_q^{-1}B P_v^{-1}A^T A_s^{-1} A P_v^{-1}B^TP_q^{-1}{\bf q}, {\bf q} \rangle}
    {4 \rho \omega \langle P_q^{-1}BP_v^{-1}B^TP_q^{-1}{\bf q}, {\bf q} \rangle}
    = \frac{\langle A_s^{-1} A {\bf w}, A {\bf w} \rangle}
    {4 \rho \omega \langle P_v{\bf w}, {\bf w} \rangle},    
  \end{equation}
  where we defined ${\bf w} = P_v^{-1}B^TP_q^{-1}{\bf q}$. Combining
  \cref{eq:Asmin1wwbound,eq:Avdualbound},
  \begin{equation}
    \label{eq:Asmin1AwAwboundn}
    \langle A_s^{-1}A{\bf w}, A{\bf w} \rangle
    \le \alpha^{-1} \langle P_v^{-1}A{\bf w}, A{\bf w} \rangle
    = \alpha^{-1} \norm[0]{A{\bf w}}_{P_v^{-1}}^2
    \le \alpha^{-1} c_a^2 \norm[0]{{\bf w}}_{P_v}^2
    = \alpha^{-1} c_a^2 \langle P_v {\bf w}, {\bf w} \rangle.
  \end{equation}
  Using \cref{eq:Asmin1AwAwboundn} in \cref{eq:lambdamin1}, we find
  that $|\lambda - 1| \le \frac{c_a}{2\sqrt{\rho \alpha
      \omega}}$. Combined with \cref{eq:infcstar}, the expression from
  $c_{\star}$ follows. Positivy of $c_{\star}$ is a consequence of
  $\rho \omega > \tfrac{1}{4}c_a^2\alpha^{-1}$.
  \\
  {\bf Step 3.} We determine $\bar{c}_{\star}$. By definition and the
  substitution ${\bf x} = \mathcal{A}{\bf y}$
  \begin{equation}
    \label{eq:starwithcstarbardeff}
    \bar{c}_{\star}
    \le \inf_{{\bf 0} \ne {\bf x} \in \mathbb{R}^n} \frac{\langle \mathcal{P}^{-1}\mathcal{P}_UA^{-1}{\bf x}, {\bf x} \rangle}{\langle \mathcal{P}^{-1}{\bf x}, {\bf x} \rangle}
    = \inf_{{\bf 0} \ne {\bf y} \in \mathbb{R}^n} \frac{\langle \mathcal{P}^{-1}P_U{\bf y}, \mathcal{A}{\bf y} \rangle}{\langle \mathcal{P}^{-1}\mathcal{A}{\bf y}, \mathcal{A}{\bf y} \rangle}.
  \end{equation}
  We consider separately the numerator and denominator. For
  ${\bf y} = [{\bf v}^T\ {\bf q}^T]^T$ and defining
  ${\bf w} = P_v^{-1}B^T{\bf q}$, we find for the numerator
  \begin{equation}
    \label{eq:numeratoreq}
    \begin{split}
      \langle \mathcal{P}^{-1} P_U{\bf y}, \mathcal{A}{\bf y} \rangle    
      =&\langle {\bf v} + \omega^{-1}P_v^{-1}B^T{\bf q}, A{\bf v}+B^T{\bf q} \rangle
         - \langle {\bf q}, B{\bf v} \rangle
      \\
      =&\langle A_s{\bf v}, {\bf v} \rangle
         + \omega^{-1}\langle A{\bf v}, {\bf w} \rangle
         + \omega^{-1}\langle {\bf w}, P_v{\bf w} \rangle.      
    \end{split}
  \end{equation}
  By Cauchy--Schwarz, Young's inequality, and \cref{eq:Avdualbound},
  \begin{equation}
    \label{eq:Avvlower}
    \langle A{\bf v}, {\bf w}\rangle \ge - \tfrac{1}{2}
    c_a^2\norm[0]{{\bf v}}_{P_v}^2 - \tfrac{1}{2} \norm[0]{{\bf
        w}}_{P_v}^2. 
  \end{equation}
  Using \cref{eq:abs-coercivity} and \cref{eq:Avvlower} to bound the
  right hand side of \cref{eq:numeratoreq} from below,
  \begin{equation}
    \label{eq:expandnumeratorff}
    \begin{split}
      \langle \mathcal{P}^{-1} P_U{\bf y}, \mathcal{A}{\bf y} \rangle
      &\ge \alpha \norm[0]{{\bf v}}_{P_v}^2
        - \tfrac{1}{2} c_a^2 \omega^{-1} \norm[0]{{\bf v}}_{P_v}^2
        - \tfrac{1}{2} \omega^{-1} \norm[0]{{\bf w}}_{P_v}^2
        + \omega^{-1}\norm[0]{{\bf w}}_{P_v}^2
      \\
      &= \del[1]{\alpha - \tfrac{1}{2} c_a^2 \omega^{-1}} \norm[0]{{\bf v}}_{P_v}^2
        + \tfrac{1}{2} \omega^{-1} \norm[0]{{\bf w}}_{P_v}^2.      
    \end{split}
  \end{equation}
  Expanding now the denominator of \cref{eq:starwithcstarbardeff},
  \begin{equation}
    \label{eq:expanddenomPAy}
    \langle \mathcal{P}^{-1}\mathcal{A}{\bf y}, \mathcal{A}{\bf y} \rangle
    = \omega^{-1}\norm[0]{A{\bf v} + B^T{\bf q}}_{P_v^{-1}}^2 + \rho^{-1}\norm[0]{B{\bf v}}_{P_q^{-1}}^2.
  \end{equation}
  Define ${\bf w}$ such that $B^T{\bf q} = P_v{\bf w}$. Using the
  triangle inequality, \cref{eq:Avdualbound,eq:Bvdualbound}, we bound
  the right hand side of \cref{eq:expanddenomPAy},
  \begin{equation}
    \label{eq:expanddenomff}
    \begin{split}
      \langle \mathcal{P}^{-1}\mathcal{A}{\bf y}, \mathcal{A}{\bf y} \rangle
      &= \omega^{-1}\norm[0]{A{\bf v} + P_v{\bf w}}_{P_v^{-1}}^2 + \rho^{-1}\norm[0]{B{\bf v}}_{P_q^{-1}}^2
      \\
      &\le 2\omega^{-1} c_a^2 \norm[0]{{\bf v}}_{P_v}^2
        + 2\omega^{-1}\norm[0]{{\bf w}}_{P_v}^2
        + \rho^{-1} c_b^2\norm[0]{{\bf v}}_{P_v}^2.      
    \end{split}
  \end{equation}
  Combining \cref{eq:expandnumeratorff,eq:expanddenomff} with
  \cref{eq:starwithcstarbardeff},
  \begin{equation*}
      \inf_{{\bf 0} \ne {\bf y} \in \mathbb{R}^n} \frac{\langle \mathcal{P}^{-1}P_U{\bf y}, \mathcal{A}{\bf y} \rangle}{\langle \mathcal{P}^{-1}\mathcal{A}{\bf y}, \mathcal{A}{\bf y} \rangle}
      \ge \inf_{{\bf 0} \ne {\bf y} \in \mathbb{R}^n} \frac{\del[1]{2\alpha\omega - c_a^2} \norm[0]{{\bf v}}_{P_v}^2 + \norm[0]{{\bf w}}_{P_v}^2}
        {(4 c_a^2+ 2\omega \rho^{-1} c_b^2) \norm[0]{{\bf v}}_{P_v}^2
        + 4\norm[0]{{\bf w}}_{P_v}^2} \ge \bar{c}_{\star}.
  \end{equation*}
  Positivity of $\bar{c}_{\star}$ is a consequence of
  $\omega > \tfrac{1}{2}c_a^2 \alpha^{-1}$.  
\end{proof}

\section{An equality used in the proof of
  \cref{lem:charac-st-preconditioner}}

\begin{lemma}
  \label{lem:inf-equality}
  The following equality holds:
  \begin{multline}
    \label{eq:inf-equality}
    \int_{I_n} \inf_{p_h = p_{h,1} + p_{h,0}}
    \del[1]{\tau^{-1} (\nabla p_{h,1}, \nabla p_{h,1})
      + \nu^{-1} (p_{h,0}, p_{h,0})
    } \kappa_n \dif{t}\\
    = 
    \inf_{p_h = p_{h,1} + p_{h,0}}
    \del[1]{
      \tau^{-1} \int_{I_n} (\nabla p_{h,1}, \nabla p_{h,1 }) \kappa_n \dif{t}
      + \nu^{-1} \int_{I_n} (p_{h,0}, p_{h,0}) \kappa_n \dif{t}
    }.
  \end{multline}
\end{lemma}
\begin{proof}
  \textbf{Step 1.} Consider the left hand side of
  \cref{eq:inf-equality}.  Let $\cbr[0]{\phi_i}_{i=0}^{k_t}$ be an
  orthonormal basis of $P_k(I_n)$ with respect to the weighted
  $L^2(I_n)$ inner product. Then
  $p_h(x,t) = \sum_{i=0}^{k_t} \phi_i(t) p_i(x)$ and we write
  \begin{align*}
    &
      \int_{I_n} \inf_{p_h = p_{h,1} + p_{h,0}}
      \del[1]{\tau^{-1} (\nabla p_{h,1}, \nabla p_{h,1})
      + \nu^{-1} (p_{h,0}, p_{h,0})
      } \kappa_n \dif{t}
    \\
    =& \sum_{i=0}^{k_t}\inf_{p_i = p_{i,1} + p_{i,0}}
       \del[2]{\tau^{-1} (\nabla p_{i,1}, \nabla p_{i,1})
       + \nu^{-1} (p_{i,0}, p_{i,0})
       }.
  \end{align*}
  The infimum for each $i$ is attained by the unique
  $p_{i,1}^* \in Q_h$ that satisfies \cite[Section
  5.4]{henriquez2026robust}
  \begin{equation*}
    \tau^{-1}(\nabla p_{i,1}^*,\nabla r_i) + \nu^{-1}(p_{i,1}^*, r_i)
    = \nu^{-1}(p_i, r_i),
  \end{equation*}
  for all $r_i \in Q_h$. Summing over $i$, and using that
  $p_{h,1}^*(x,t) = \sum_{i=0}^{k_t} \phi_i(t) p_{i,1}^*(x) \in
  \mathcal{Q}_h$ and
  $r_h(x,t) = \sum_{i=0}^{k_t} \phi_i(t) r_i(x) \in \mathcal{Q}_h$,
  \begin{equation}
    \label{eq:charph1star}
    \int_{I_n}\del[2]{\tau^{-1}(\nabla p_{h,1}^*,\nabla r_h) + \nu^{-1}(p_{h,1}^*, r_h)} \kappa_n \dif t
    = \int_{I_n} \nu^{-1}(p_h, r_h) \kappa_n \dif t.
  \end{equation}
  \textbf{Step 2.} The infimum of the right hand side of
  \cref{eq:inf-equality} is attained by the unique
  $p_{h,1}^{\dagger} \in \mathcal{Q}_h$ that satisfies
  \begin{equation}
    \label{eq:charph1dag}
    \tau^{-1} \int_{I_n} (\nabla p_{h,1}^{\dagger}, \nabla s_h) \kappa_n \dif{t}
    + \nu^{-1} \int_{I_n} (p_{h,1}^{\dagger}, s_h) \kappa_n \dif{t}
    = \nu^{-1} \int_{I_n} (p_{h}, s_h) \kappa_n \dif{t},
  \end{equation}
  for all $s_h \in \mathcal{Q}_h$. Since
  \cref{eq:charph1star,eq:charph1dag} are the same equation, we
  conclude that $p_{h,1}^* = p_{h,1}^{\dagger}$ and so
  \cref{eq:inf-equality} holds.
\end{proof}

\end{document}